\documentclass[11pt,reqno]{preprint}

\usepackage[full]{textcomp}
\usepackage[osf]{newtxtext}
\usepackage{comment}

\usepackage{amssymb}
\usepackage{mathtools}
\usepackage{hyperref}
\usepackage{breakurl}
\usepackage{mhenvs}
\usepackage{mhequ}
\usepackage{mhsymb}
\usepackage{booktabs}
\usepackage{tikz}
\usetikzlibrary{decorations.pathreplacing}
\usetikzlibrary{decorations.markings}
\usetikzlibrary{decorations.pathmorphing}
\usetikzlibrary{snakes}
\usepackage{mathrsfs}
\usepackage{longtable}
\usepackage{scalerel}
\usepackage{cprotect}
\usepackage{multirow}
\usepackage{colortbl}
\usepackage{makecell}
\usepackage{resizegather}
\usepackage{enumerate}

\definecolor{ink}{HTML}{203040}
\definecolor{blue}{HTML}{2563EB}
\definecolor{teal}{HTML}{0F8B8D}
\definecolor{orange}{HTML}{E79A36}
\definecolor{red}{HTML}{C44D58}
\definecolor{pale}{HTML}{F2F5F8}
\definecolor{muted}{HTML}{667085}
\definecolor{outsidefill}{HTML}{F0F4F8}
\definecolor{insidefill}{HTML}{DDF2EB}
\definecolor{otherfill}{HTML}{FFF0D9}
\definecolor{grid}{HTML}{CBD4DD}
\definecolor{wall}{HTML}{D28B30}
\definecolor{cycle}{HTML}{BE4350}

\usetikzlibrary{arrows.meta,calc}

\newcommand{\SU}{\mathrm{SU}}
\newcommand{\dH}{\mathrm dH}
\newcommand{\dd}{\mathrm d}

\newcommand{\Prob}{\mathbb P}

\def\bt{\begin{theorem}}
	\def\et{\end{theorem}}
\def\bl{\begin{lemma}}
	\def\el{\end{lemma}}
\def\br{\begin{remark}}
	\def\er{\end{remark}}
\def\bx{\begin{Examples}}
	\def\ex{\end{Examples}}
\def\bd{\begin{definition}}
	\def\ed{\end{definition}}
\def\bp{\begin{proposition}}
	\def\ep{\end{proposition}}
\def\bc{\begin{corollary}}
	\def\ec{\end{corollary}}
\def\<{{\langle}}
\def\>{{\rangle}}

\newcommand{\Tr}{\operatorname{Tr}}
\def\Re{\operatorname{Re}}

\newcommand{\Cov}{\operatorname{Cov}}

\definecolor{codegreen}{rgb}{0,0.6,0}
\definecolor{codegray}{rgb}{0.5,0.5,0.5}
\definecolor{codepurple}{rgb}{0.58,0,0.82}
\definecolor{backcolour}{rgb}{0.95,0.95,0.92}

\tikzset{
    midarrow/.style={postaction={decorate,decoration={markings,mark=at position 0.7 with {\arrow{>}}}}},
    midarrow1/.style={postaction={decorate,decoration={markings,mark=at position 0.45 with {\arrow{>}}}}},
    midarrow2/.style={postaction={decorate,decoration={markings,mark=at position 0.75 with {\arrow{>}}}}},
}

\colorlet{darkblue}{blue!90!black}
\colorlet{darkred}{red!90!black}
\colorlet{darkgreen}{green!70!black}

\title{Anisotropic lattice Yang--Mills:
mass gap   at semi-strong couplings}

\author{Hao Shen$^1$}

\institute{University of Chicago, \email{pkushenhao@gmail.com}}

\begin{document}

\maketitle

\begin{abstract}
The mass gap problem for the lattice Yang--Mills model at any inverse couplings $\beta$ is so far an open problem.
For dimensions three or higher, it was proved only in the strong coupling regime, namely for $\beta\ge 0$ small.
In this paper we introduce the anisotropic lattice Yang--Mills model,
with a parameter $\beta_{ij}\ge 0$ for each plaquette direction $i<j$.
We first show that in 3D, mass gap holds for any $\beta_{12}$,
and small $\beta_{13},\beta_{23}$ whose smallness condition does not depend on the value of $\beta_{12}$.
The proof easily extends to higher dimensions, when $\beta_{12}$, $\beta_{34}$, etc. are arbitrary and the other $\beta_{ij}$ are small. 
We then show uniqueness of the infinite-volume Gibbs measure.
Finally, in 3D, we discuss a special case where $\beta_{12}$ is arbitrary and ``half'' of the plaquette-dependent $\beta_{13},\beta_{23}$ are arbitrary, in an alternating or staggered pattern,  for which we can also prove the mass gap.
\end{abstract}

\setcounter{tocdepth}{2}
\tableofcontents

\section{Introduction and main result}

Let $N\ge2$, let $G=\SU(N)$, and let $H$ be Haar probability measure on $G$.
We will work in a finite rectangular box $\Lambda\subset\Z^n$ with free 
boundary conditions.
Associate $Q_e\in G$ to each positively oriented edge $e\in E^+_\Lambda$ and set
$Q_{e^{-1}}=Q_e^{-1}$. 
For each unit square choose one starting vertex and one orientation,
using the lexicographic convention; denote the resulting
set of plaquettes by $P^+_\Lambda$. 

For a plaquette $p=e_1e_2e_3e_4$, write
\[
 Q_p=Q_{e_1}Q_{e_2}Q_{e_3}Q_{e_4},\qquad W_p=\Tr Q_p,\qquad w_p=N^{-1}W_p ,
\]
where $\Tr$ is the usual matrix trace. A plaquette $p$ parallel to the $ij$-plane has
coupling $\beta_p=\beta_{ij}\ge 0$. The Wilson Yang--Mills measure is
\begin{equation}\label{eq:measure}
\dd\mu(Q)=
 \dd\mu_{\Lambda,N,\boldsymbol\beta}(Q)
 =Z_{\Lambda,N,\boldsymbol\beta}^{-1}
 \exp\!\Big(
 N\sum_{p\in P^+_\Lambda}\beta_p\operatorname{Re}\Tr Q_p
 \Big)
 \prod_{e\in E^+_\Lambda}\dH(Q_e).
\end{equation}
For complex observables $F,G$, write
\[
 \Cov_\mu(F,G)=\mu(F\overline G)-\mu(F)\,\overline{\mu(G)} .
\]
Let $d(p,p')$ be the minimum $\ell^1$ distance between a vertex of $p$ and a vertex
of $p'$.

For clearer exposition, we first state and prove our result in 3D. Then we extend the result to higher dimensions.

\begin{theorem}\label{theo:main}
Let $n=3$.
There are absolute constants $b_*,C>0$ 
such that for every $N\ge2$, all couplings satisfying
\begin{equation}\label{eq:small}
 \beta_{12}<\infty,\qquad
  N^2\max(\beta_{13},\beta_{23})\le b_* ,
\end{equation}
and any two plaquettes $p,p'$ in a free-boundary rectangular box $\Lambda$,  one has
\begin{equation}\label{eq:main}
 |\Cov_\mu(W_p,W_{p'})|\le CN^2 e^{-m d(p,p')},\qquad m=\tfrac12\log 2.
\end{equation}
The same estimate holds if one trace is conjugated or if both traces are replaced by
their real parts. 
The constants are independent of $N$, $\beta_{12}$ and the size of $\Lambda$.
%
\end{theorem}

The proof of Theorem~\ref{theo:main} is given in the end of Section~\ref{sec:counting}.
For standard (i.e. isotropic) lattice Yang--Mills models where $\beta_{ij}=\beta$ is sufficiently small, this can be proved by cluster expansion \cite{Seiler,OS1978};
exponential decay was also proved in \cite{SZZ22,Nissim25} 
assuming that $\beta<\beta_0$ for explicit $\beta_0$ uniform in $N$.

Now we extend the result to
higher dimensions.
Let  $n\ge 4$. After permuting the coordinates, take the
$\lfloor n/2\rfloor$ {\it pairs of orientations} 
\[
(2a-1,2a),\qquad 
1\le a\le\lfloor n/2\rfloor.
\]
For example,
in dimension four,  the  couplings in the paired orientations are $\beta_{12},\beta_{34}$,
and the other orientations  are $13,14,23,24$.

\begin{theorem}\label{theo:high-dim}
Let $n\ge 4$ be the lattice dimension. 
There are constants $b_*(n),C_n>0$ such that, for every $N\ge2$,
if the couplings 
\[
 \beta_{2a-1,2a} <\infty,\quad (1\le a\le\lfloor n/2\rfloor), \qquad
 N^2\beta_r\le b_*(n)
\]
 for every $r$ not in the paired orientations, then
\begin{equation}\label{eq:disjoint-pairs-decay}
  \big|\operatorname{Cov}_{\mu}(W_p,W_{p'})\big|
  \le C_n N^2 e^{-m d(p,p')},\qquad m=\tfrac12\log 2 .
\end{equation}
The estimate is uniform in $\Lambda$ and all the couplings, and holds also
with one trace conjugated or with both traces replaced by their real parts.
\end{theorem}

The proof of Theorem~\ref{theo:high-dim} is given in Section~\ref{sec:Higher-dim} and is essentially the same as that of Theorem~\ref{theo:main}. We choose to first demonstrate the proof in 3D because it is easier to visualize.

Once we show exponential decay, with minor modification of the proof plus essentially standard DLR argument, we can prove the Gibbs uniqueness of infinite volume limit. This is given in Theorem~\ref{theo:gibbs-uniqueness}.

Finally, we discuss another extension in 3D, which allows ``half''  of the plaquette dependent inverse couplings in the 13 and 23 directions to be arbitrary, in an ``alternating'' or ``staggered'' pattern.

Let $n=3$.
At a plaquette
with minimum vertex $(i,j,k)$, set
\begin{equation}\label{eq:ym3-ring-pattern}
 \beta_{12}=\beta_A,\qquad
 \beta_{13}=
 \begin{cases}
   \beta_B & i\text{ even},\\
   \beta_b & i\text{ odd},
 \end{cases}
 \qquad
 \beta_{23}=
 \begin{cases}
   \beta_C & j+k\text{ even},\\
   \beta_c & j+k\text{ odd}.
 \end{cases}
\end{equation}
All five couplings are finite and nonnegative. An uppercase subscript
denotes a coupling that can be arbitrarily large. Each vertical orientation
has half of its plaquettes in the unrestricted class in one period cell.
The $13$ pattern repeats at every height. The $23$ pattern shifts by one
row at each step in direction~$3$.

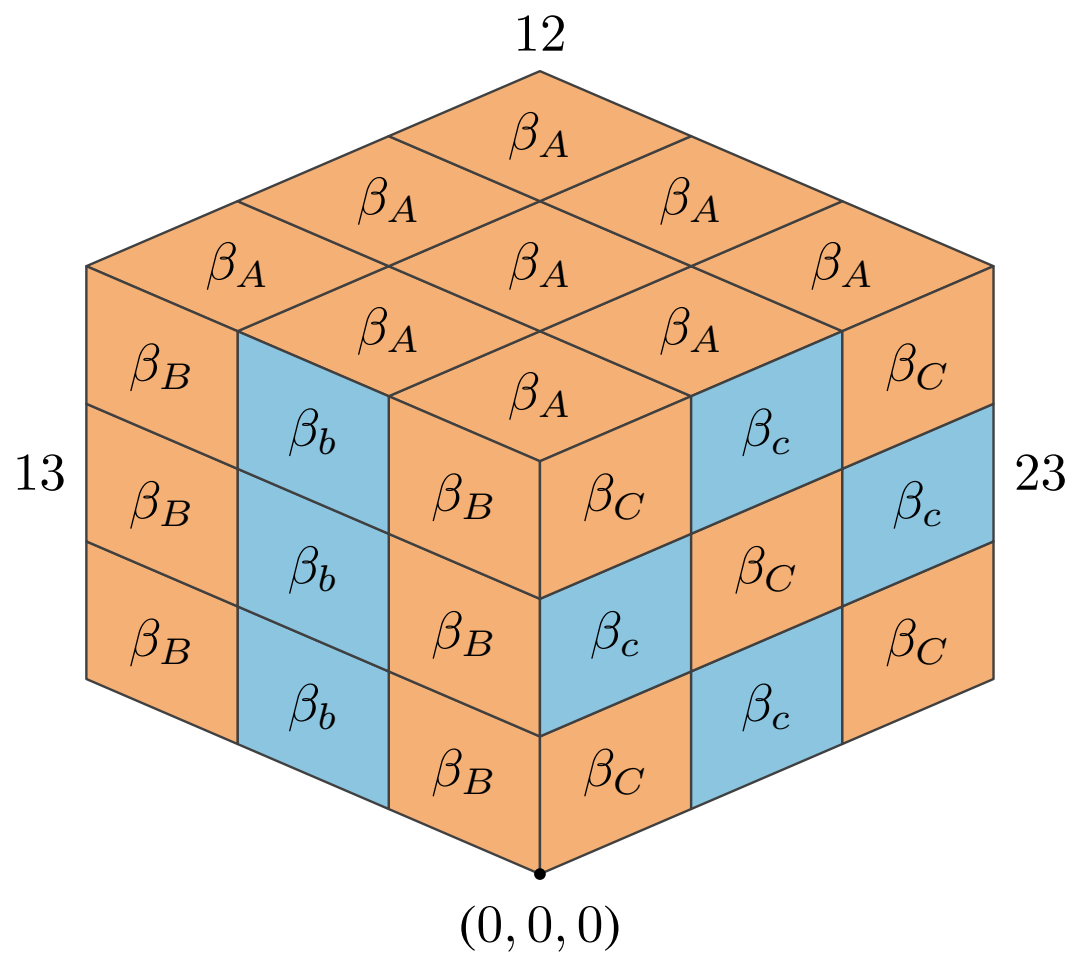
\begin{figure}[htbp]
\centering
\begin{tikzpicture}[
  x={(-1cm,.43cm)}, y={(1cm,.43cm)}, z={(0cm,.91cm)},
  tile/.style={draw=black!75, line width=.45pt, line join=round},
  font=\small
]
  \definecolor{ymOrange}{RGB}{244,176,117}
  \definecolor{ymBlue}{RGB}{139,197,224}

  \foreach \u in {0,1,2} {
    \foreach \v in {0,1,2} {
      \path[tile,fill=ymOrange] (\u,\v,3) -- ++(1,0,0)
        -- ++(0,1,0) -- ++(-1,0,0) -- cycle;
      \node at ({\u+.5},{\v+.5},3) {$\beta_A$};

      \ifodd\u\relax
        \def\ymFill{ymBlue}\def\ymIndex{b}
      \else
        \def\ymFill{ymOrange}\def\ymIndex{B}
      \fi
      \path[tile,fill=\ymFill] (\u,0,\v) -- ++(1,0,0)
        -- ++(0,0,1) -- ++(-1,0,0) -- cycle;
      \node at ({\u+.5},0,{\v+.5}) {$\beta_{\ymIndex}$};

      \ifodd\numexpr\u+\v\relax
        \def\ymFill{ymBlue}\def\ymIndex{c}
      \else
        \def\ymFill{ymOrange}\def\ymIndex{C}
      \fi
      \path[tile,fill=\ymFill] (0,\u,\v) -- ++(0,1,0)
        -- ++(0,0,1) -- ++(0,-1,0) -- cycle;
      \node at (0,{\u+.5},{\v+.5}) {$\beta_{\ymIndex}$};
    }
  }
  \node[above] at (3,3,3) {$12$};
  \node[left]  at (3,0,1.5) {$13$};
  \node[right] at (0,3,1.5) {$23$};
  \fill (0,0,0) circle[radius=1.1pt];
  \node[below=2pt] at (0,0,0) {$ (0,0,0) $};
\end{tikzpicture}
\caption{Staggered couplings on a $3\times3\times3$ cube.
Orange plaquettes carry arbitrary $\beta_A,\beta_B,\beta_C$.
Blue plaquettes carry small  $\beta_b,\beta_c$.
}
\label{fig:ym3-staggered-cube}
\end{figure}

\begin{theorem}\label{theo:3D-ring}
Let $n=3$.
There are absolute constants $\delta,C,m>0$ such that
for every $N\ge2$ and
\[
\beta_A,\beta_B,\beta_C <\infty,
\qquad
N^2\max\{\beta_b,\beta_c\}\le\delta
\]
we have
\begin{equation}\label{eq:ym3-ring-decay}
 \big|\operatorname{Cov}_{\mu}(W_p,W_{p'})\big|
 \le C N^2 e^{-m d(p,p')}
\end{equation}
for every pair of plaquettes $p,p'$ in $\Lambda$. 
The same bound holds if one trace is conjugated
or if both are replaced by their real parts. 
Infinite volume sequences give  the
same local limit and \eqref{eq:ym3-ring-decay} holds in that limit.
\end{theorem}

The proof of Theorem~\ref{theo:3D-ring} is given in Section~\ref{sec:ym3-ring}, which is the same method as in the proof of Theorem~\ref{theo:main} with an additional trick.  

Our proof of the above theorems relies on solvability structure in 2D Yang--Mills (c.f. \cite{Driver89,Levy11}), and a replica and switching idea used earlier for Ising model and percolation, c.f. 
\cite{GriffithsHurstSherman1970,Percus1975,Aizenman1982,vanDenBergMaes1994}, or \cite[Section~3]{DuminilCopin2018},
and also used in lattice Yang--Mills with finite groups \cite{For22,AC25}.
The switching idea reduces correlation problems to connecting probabilities, often achieved by some type of ``flipping'' (though not necessarily by swapping replica) in  the component of one observable, c.f. \cite{EdwardsSokal1988,ChayesMachtaRedner1998,CampbellChayes1998}. We make some comments on our results.
\begin{enumerate}
\item
The delicate design of the staggered pattern in Theorem~\ref{theo:3D-ring} is only to demonstrate that in some special cases, we are able to not just integrate out edge variables in the horizontal layers using spanning tree argument (as in Theorem~\ref{theo:main}) but also integrate out some of the vertical edges as long as this only creates ``relatively small and separated'' holonomies.
It is very likely that such ``staggered results'' can be extended to higher dimensions, but it seems that such extensions are not very essential, and somewhat artificial, so we do not study this further.
\item
Remark that the smallness condition in \eqref{eq:small} depends on $N$.
In \cite{SZZ22}, the small $\beta$ condition is independent of $N$, which allows to study large $N$ problems in that setting. Here however our motivation is not large $N$ problems; the truly interesting question here is that even for the weak and strong forces where $N=2$ and $N=3$ respectively, what is the minimum number of $\beta_{ij}$ that has to be assumed small for the proof to go through.
\item
To simplify notation, we only study the covariance of $W_p,W_{p'}$ here. It is completely straightforward to adapt the proof to the case of local bounded observables $F,G$. 
This is essentially treated in the proof of Theorem~\ref{theo:gibbs-uniqueness} for one local bounded observable $F$.
In particular our proof does not require that the observables are gauge invariant.
\item
Our method does not seem to apply to the case that all $\beta_p$ are large. It is not clear yet if the method can be used to prove the result with $\beta_{12}$, $\beta_{23}$ arbitrary  and only $\beta_{13}$ small in 3D.
\end{enumerate}

Lattice gauge theories with anisotropic or plaquette-dependent couplings appeared in physics literature, c.f. \cite{karsch1982n,billo1997,Orland07,giordano2015}, but there has not been mathematically rigorous study as far as we know.
(The paper \cite{Nissim25} mentioned that $\beta$ is allowed to depend on plaquettes $p$, which is convenient  in his proof of infinite volume limit, though the result therein only assumes $\sup_p \beta$ bounded by a small constant.)

\medskip

{\bf Acknowledgements:} H.S. gratefully acknowledges financial support from the NSF grant CAREER
DMS-2044415 / DMS-2706725.

\medskip

{\bf AI statement:} ChatGPT 6.2 and Gemini were used to draw all the pictures, help brainstorm the staggered pattern,  check the counting arguments, and polish the paper in the final round. 

\section{Setup for proof}
\label{sec:setup}

Recall that $\mu$ is the Yang--Mills measure \eqref{eq:measure}. A replica is a full
configuration of edge variables sampled from $\mu$. Draw $Q^1$ and $Q^2$
independently, each with this law. Their joint law is $\mu\otimes\mu$.

For two  observables $F$ and $G$, write $F^a=F(Q^a)$ and $G^a=G(Q^a)$
for $a=1,2$. In the application, $F=w_p$ and $G=w_{p'}$. The covariance
can be written as
\begin{equation}\label{eq:rep}
 2\Cov_\mu(F,G)=
 \E_{\mu\otimes\mu}\!\left[(F^1-F^2)\overline{(G^1-G^2)}\right].
\end{equation}
Indeed this holds by expanding the product inside the expectation.

We call the $12$-plaquettes {\it horizontal}. 
Let $T$ be the set of all the $13$- and
$23$-plaquettes, which we call {\it vertical}. 
The following picture shows a vertical plaquette $r\in T$.
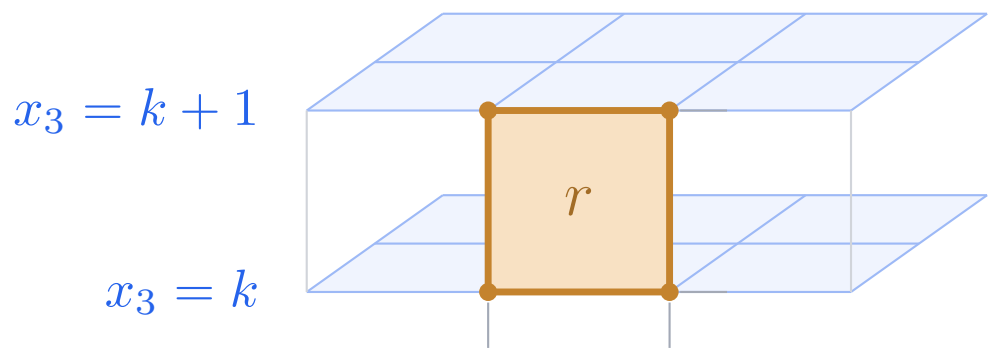
\[
\begin{tikzpicture}[x=1cm,y=1cm,font=\small]
  \begin{scope}[shift={(.2,.85)},x={(1.2cm,0cm)},y={(.45cm,.32cm)},z={(0cm,1.2cm)}]
    \foreach \z in {0,1}{
      \fill[blue!7] (0,0,\z)--(3,0,\z)--(3,2,\z)--(0,2,\z)--cycle;
      \foreach \a in {0,1,2,3}{\draw[blue!45] (\a,0,\z)--(\a,2,\z);}
      \foreach \a in {0,1,2}{\draw[blue!45] (0,\a,\z)--(3,\a,\z);}
    }
    \draw[muted!30] (0,0,0)--(0,0,1) (3,0,0)--(3,0,1);
    \filldraw[fill=orange!30,draw=orange!85!black,line width=1.3pt]
      (1,0,0)--(2,0,0)--(2,0,1)--(1,0,1)--cycle;
    \foreach \x/\z in {1/0,2/0,2/1,1/1}{\fill[orange!85!black] (\x,0,\z) circle (1.7pt);}
    \node[orange!70!black,font=\large] at (1.5,0,.5) {$r$};
  \end{scope}
  \node[blue,anchor=east,font=\small] at (.02,2.05) {$x_3=k+1$};
  \node[blue,anchor=east,font=\small] at (.02,.85) {$x_3=k$};
  \draw[muted!60] (1.4,.78)--(1.4,.48) (2.6,.78)--(2.6,.48);
  \draw[muted!60] (2.67,.85)--(2.98,.85) (2.67,2.05)--(2.98,2.05);
\end{tikzpicture}
\]
Now we describe a process to ``mark'' the vertical plaquettes, as follows.
Only the plaquettes in $T$ will be
marked. For any plaquette $f$, write
\[
 t_f=N^2\beta_f \ge 0,\qquad 
 X_f^a=\operatorname{Re}w_f(Q^a) \in [-1,1]
\]
since $|w_f(Q^a)|\le1$
because $Q_f^a$ is unitary.
A vertical plaquette $r$ contributes the factor $e^{t_r X_r^a}$ to the
weight. Write 
\begin{equ}[e:expand-H]
 e^{t_r(X_r^1+X_r^2)}=e^{-2t_r}(1+H_r),\qquad
 H_r=e^{t_r(2+X_r^1+X_r^2)}-1 \ge 0.
\end{equ}
We can now interpret $1$ and $H_r$ as the two
possible weights of a mark.

Write $\mathbf Q=(Q^1,Q^2)$. Once $\mathbf Q$ has been sampled, each
$H_r=H_r(\mathbf Q)$ is a known number. 
At each $r\in T$, independently
toss a coin whose two probabilities are
\[
 \Prob(r\text{ is marked}\mid\mathbf Q) =\frac{H_r}{1+H_r},\qquad
 \Prob(r\text{ is unmarked}\mid\mathbf Q)=\frac{1}{1+H_r}.
\]
Let $R$ be the set of \emph{marked}  plaquettes, and let $\Prob$
denote the joint law of 
\[
(Q^1,Q^2,R).
\]
Note that independence of the coins is only conditional on $\mathbf Q$.

For  $S\subseteq T$, we have
\begin{equ}[e:PSQ]
 \Prob(R=S\mid\mathbf Q)
 =\prod_{r\in S}\frac{H_r}{1+H_r}
     \prod_{r\in T\setminus S}\frac{1}{1+H_r}
 =\frac{\prod_{r\in S}H_r}{\prod_{r\in T}(1+H_r)}.
\end{equ}
As a sanity check, the sum over all $S\subseteq T$ of \eqref{e:PSQ} is one, and for any event $A$
depending only on the two fields,
\[
 \Prob(\mathbf Q\in A)
 =\int_A\bigg(\sum_{S\subseteq T}\Prob(R=S\mid\mathbf Q)\bigg)
           (\mu\otimes\mu)(\dd\mathbf Q)
 =(\mu\otimes\mu)(A).
\]
Thus the marginal law of the two fields is the product law $\mu\otimes\mu$ used in
\eqref{eq:rep}.

Again conditioned on $\mathbf Q$, 
for a fixed $S$, the probability that $R$ contains $S$ is
\[
 \Prob(R\supseteq S\mid\mathbf Q)=
 \mathbb E \Big[
 \prod_{r\in S}1_{r\in R}\mid\mathbf Q\Big]
 =\prod_{r\in S}\frac{H_r}{1+H_r}.
\]
By definition we have
\[
 \frac{H_r}{1+H_r}=1-e^{-t_r(2+X_r^1+X_r^2)}\le 1-e^{-4t_r}.
\]
The upper bound holds
for every pair of fields $\mathbf Q$, so averaging over $\mathbf Q$ yields
\begin{equation}\label{eq:dom}
 q:=1-\exp\!\left[-4N^2\max_{r\in T}\beta_r\right],\qquad
 \Prob(R\supseteq S)\le q^{|S|}\quad(S\subseteq T).
\end{equation}
Remark that this estimate also allows
plaquette-dependent couplings with the same upper bound on the vertical plaquettes. In particular, $\beta_{12}$
does not enter $q$.

We also
need the reverse conditional law, that of the fields given $R=S$.
Write 
\[
\nu(\dd\mathbf Q)=\prod_{e\in E^+_\Lambda}\dH(Q_e^1)\dH(Q_e^2)
\]
for product Haar measure on the two configurations, and write the
product of the horizontal weights
\[
 L(\mathbf Q)=\prod_{f\ \mathrm{horizontal}} e^{t_f(X_f^1+X_f^2)}.
\]
We have, by \eqref{e:expand-H}
\[
\frac{\dd (\mu\otimes\mu)}{\dd \nu}(\mathbf Q)
=
Z_{\Lambda,N,\boldsymbol\beta}^{-2}
\prod_{f\notin T}
e^{t_f(X_f^1+X_f^2)}
\prod_{f\in T}
e^{t_f(X_f^1+X_f^2)}
=
c_{\Lambda,N,\boldsymbol\beta}
 L(\mathbf Q)
 \prod_{r\in T} (1+H_r)
\]
where $c_\Lambda=c_{\Lambda,N,\boldsymbol\beta}$
collects the normalizing constants and the factors $e^{-2t_r}$ appearing in \eqref{e:expand-H}.

To obtain the joint density for the fields and the event $R=S$, we multiply the above density by $ \Prob(R=S\mid\mathbf Q)$ given in \eqref{e:PSQ}:
\[
c_\Lambda L(\mathbf Q)\prod_{r\in T}(1+H_r)
\frac{\prod_{r\in S}H_r}{\prod_{r\in T}(1+H_r)}
 =c_\Lambda L(\mathbf Q)\prod_{r\in S}H_r.
\]
This joint density has a simple form which involves only the horizontal plaquette weights in $ L(\mathbf Q)$ and the marking weights $H_r$; of course, $H_r$ implicitly depends on the vertical plaquette weight.

Integrating  gives the marginal probability 
\[
\Prob(R=S)=c_\Lambda Z_S,
\qquad
Z_S:=\int L(\mathbf Q)\prod_{r\in S}H_r(\mathbf Q)\,\nu(\dd\mathbf Q).
\]
Dividing the above two lines gives
\begin{equ}\label{eq:dens}
 \Prob(\dd\mathbf Q\mid R=S)
 =\frac{1}{Z_S}
   \bigg(\prod_{r\in S}H_r\bigg)
   \bigg(\prod_{f\ \mathrm{horizontal}}e^{t_f(X_f^1+X_f^2)}\bigg)
   \nu(\dd\mathbf Q).
\end{equ}
Of course, conditioning on $R=S$, the fields $Q^1$ and $Q^2$ need not be independent: $H_r$ defined in \eqref{e:expand-H} does not allow a factorization of the above probability into two factors only depending on $Q^1$ and on $Q^2$.

Heuristically, observe that if $R=\varnothing$,  then in \eqref{eq:dens} the horizontal layers are completely decoupled, 
so $W_p$ and $W_{p'}$  are simply decorrelated if $p,p'$ are  horizontal plaquettes in two different layers, and the same holds if $p,p'$ are horizontal plaquettes in the same layer by the standard axial gauge or equivalently spanning tree argument for 2D Yang--Mills. The proof can be viewed as a ``perturbation'' of the case $R=\varnothing$: if $\beta$ in the vertical directions are small, the $q$ in \eqref{eq:dom} is small so it is rare to have a large $R$ which couple different layers.

What we use
under this conditional law \eqref{eq:dens} is an exchange symmetry. In any one factor
of \eqref{eq:dens}, swapping the two copies on all four edges of its
plaquette 
\[
Q^1_e \leftrightarrow Q^2_e \qquad  \forall e\in f
\]
 interchanges $X_f^1$ and $X_f^2$ and leaves their sum  $X_f^1+X_f^2$ unchanged.
 Of course,  when swapping the two copies on all four edges of $f$, 
since a neighboring plaquette $f'$ of $f$  only shares one edge with $f$, the other three edges of $f'$ are not swapped. 
The rest of the proof relies on two ideas:
\begin{itemize}
\item
We identify  sets of edges
on which a swap preserves every affected factor at once. 
\item
When such a
swap exchanges $F^1$ with $F^2$ while leaving $G^1,G^2$ fixed, the
integrand in \eqref{eq:rep} changes sign. This is the key to prove the decay of $\Cov_\mu(F,G)$.
\end{itemize}

\section{Integration in each horizontal layer}
\label{sec:layer}

Fix $R$ and the two observed plaquettes $p,p'$ in \eqref{eq:main}. We say that  two plaquettes in
$R\cup\{p,p'\}$ are adjacent if they share a vertex. For each connected component $A$
under this adjacency, let $B(A)$ be the smallest closed coordinate box containing its vertices. 
\begin{figure}[h]
\centering
\begin{tikzpicture}[x=.68cm,y=.68cm,font=\small,scale=0.7]
  \fill[blue!6] (0,0) rectangle (7,6);
  \fill[teal!6] (3,3) rectangle (10,9);
  \draw[blue,dashed,thick] (0,0) rectangle (7,6);
  \draw[teal,dashed,thick] (3,3) rectangle (10,9);
  \foreach \x in {0,...,6}{\filldraw[fill=blue!25,draw=blue] (\x,0) rectangle ++(1,1);}
  \foreach \y in {1,...,5}{\filldraw[fill=blue!25,draw=blue] (0,\y) rectangle ++(1,1);}
  \foreach \x in {3,...,9}{\filldraw[fill=teal!25,draw=teal] (\x,8) rectangle ++(1,1);}
  \foreach \y in {3,...,7}{\filldraw[fill=teal!25,draw=teal] (9,\y) rectangle ++(1,1);}
  \node[blue] at (2.4,1.6) {$A_1\ni p$};
  \node[teal] at (7.3,7.1) {$A_2$};
  \fill[orange!9] (13,1) rectangle (16,4);
  \draw[orange,dashed,thick] (13,1) rectangle (16,4);
  \foreach \x/\y in {13/1,14/2,15/3}{\filldraw[fill=orange!30,draw=orange] (\x,\y) rectangle ++(1,1);}
  \node[orange!75!black] at (14.5,4.7) {$A_3\ni p'$};
  \draw[blue!70!teal,line width=1.6pt] (3,3) rectangle (7,6);
\end{tikzpicture}
\caption{An illustration of $R\cup\{p,p'\}$ with
three connected components $A_1,A_2,A_3$, for which 
$B(A_1)\cap B(A_2)\ne\varnothing$. All plaquettes are in the 3D space, and $B(A_i)$ are 3D boxes, but we can only draw the picture in 2D.}
\end{figure}
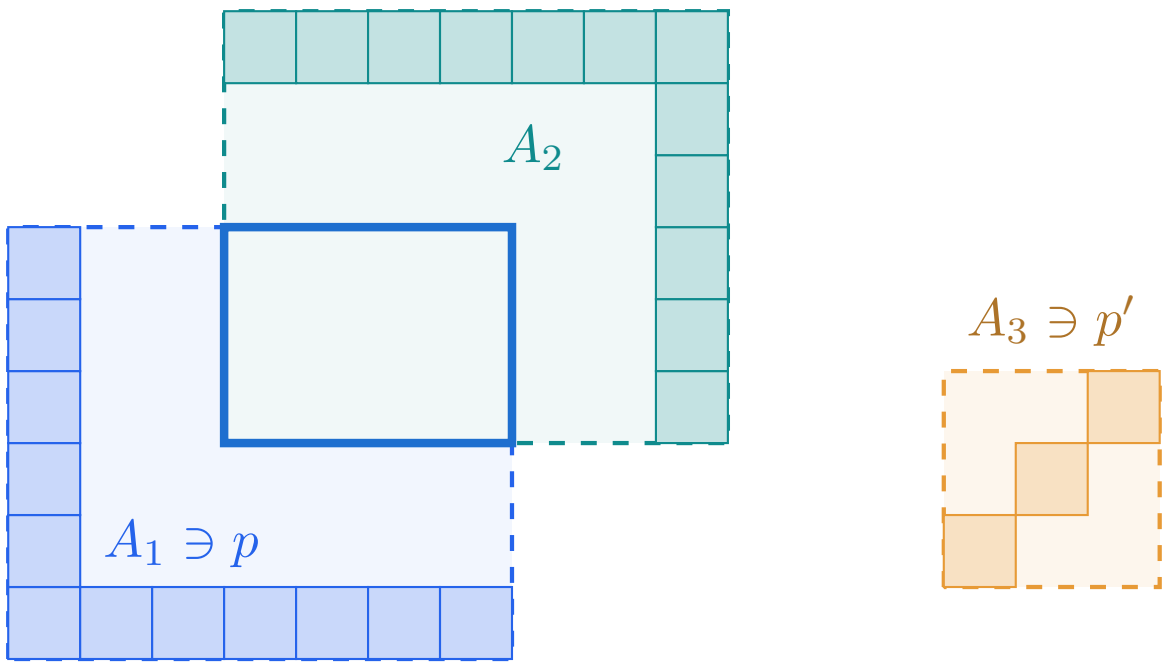

We form a graph $\mathcal G$ on these components by joining $A$ to $A'$ whenever
$B(A)\cap B(A')\ne\varnothing$. 
Write $p\leftrightarrow p'$ when the components containing $p$ and $p'$ are connected in $\mathcal G$. 

\begin{lemma}
We have
\begin{equation}\label{eq:switch}
 |\Cov(w_p,w_{p'})|\le 2\Prob(p\leftrightarrow p').
\end{equation}
\end{lemma}

\begin{proof}
Fix a horizontal layer. Let $\Gamma$ be the union of the horizontal edges used by
the plaquettes in $R\cup\{p,p'\}$ in this layer. 
We will integrate out all the edge variables $Q^1_e$, $Q^2_e$
for edges $e$ in the ``exterior of'' $\Gamma$;
 we want to keep the edge variables involved in $w_p$, $w_{p'}$ and $H_r$ untouched. To do this integration we follow an idea often used in 2D Yang--Mills which is integrating along a tree \cite{Chevyrev19YM,piazza2026log}.

{\bf Step 1.}
Consider the dual graph of the layer lattice,
with one vertex for each plaquette, and one exterior vertex $\partial$
which should be thought of as a vertex from ``very far away'';
each edge of the dual graph crosses an edge in the original lattice (sometimes called a primal edge).

\begin{figure}[h]
\centering
\begin{tikzpicture}[x=.85cm,y=.85cm,font=\small]
  \fill[teal!8] (0,0) rectangle (8,5);
  \fill[orange!17] (3,1) rectangle (5,3);
  \draw[muted!45,step=1] (0,0) grid (8,5);
  \draw[orange,very thick] (3,1) rectangle (5,3);
  \node[orange!75!black] at (4,2) {$f\notin D$};
  \node[teal] at (6.5,3.5) {$D$};
  \node[circle,draw=blue,fill=blue!9,inner sep=3pt] (root) at (-1.2,.5) {$\partial$};
  \draw[blue,thick,-{Stealth}] (root)--(.5,.5);
  \draw[blue,thick,-{Stealth}] (.5,.5)--(1.5,.5);
  \draw[blue,thick,-{Stealth}] (1.5,.5)--(2.5,.5);
  \draw[blue,thick,-{Stealth}] (1.5,.5)--(1.5,1.5);
  \draw[blue,thick,-{Stealth}] (1.5,1.5)--(1.5,2.5);
  \draw[blue,thick,-{Stealth}] (1.5,2.5)--(1.5,3.5);
  \draw[blue,thick,-{Stealth}] (1.5,3.5)--(2.5,3.5);
  \draw[blue,thick,-{Stealth}] (2.5,3.5)--(3.5,3.5);
  \draw[blue,thick,-{Stealth}] (3.5,3.5)--(4.5,3.5);
  \draw[blue,thick,-{Stealth}] (2.5,.5)--(3.5,.5);
  \draw[blue,thick,-{Stealth}] (3.5,.5)--(4.5,.5);
  \draw[red,line width=3pt] (1,0)--(1,1);
  \foreach \x/\y/\n in {.5/.5/1,1.5/.5/2,2.5/.5/3}{
    \node[circle,fill=white,draw=blue,inner sep=1.6pt,font=\scriptsize] at (\x,\y) {$\n$};
  }
\end{tikzpicture}
\caption{Picture of a horizontal layer. The orange edges are those in $\Gamma$.
The blue arrows are part of the spanning tree;
the numbers $1,2,3$ correspond to the plaquettes we process in this order; the red edge is the edge variable we integrate when we process the plaquette 2.}
\end{figure}
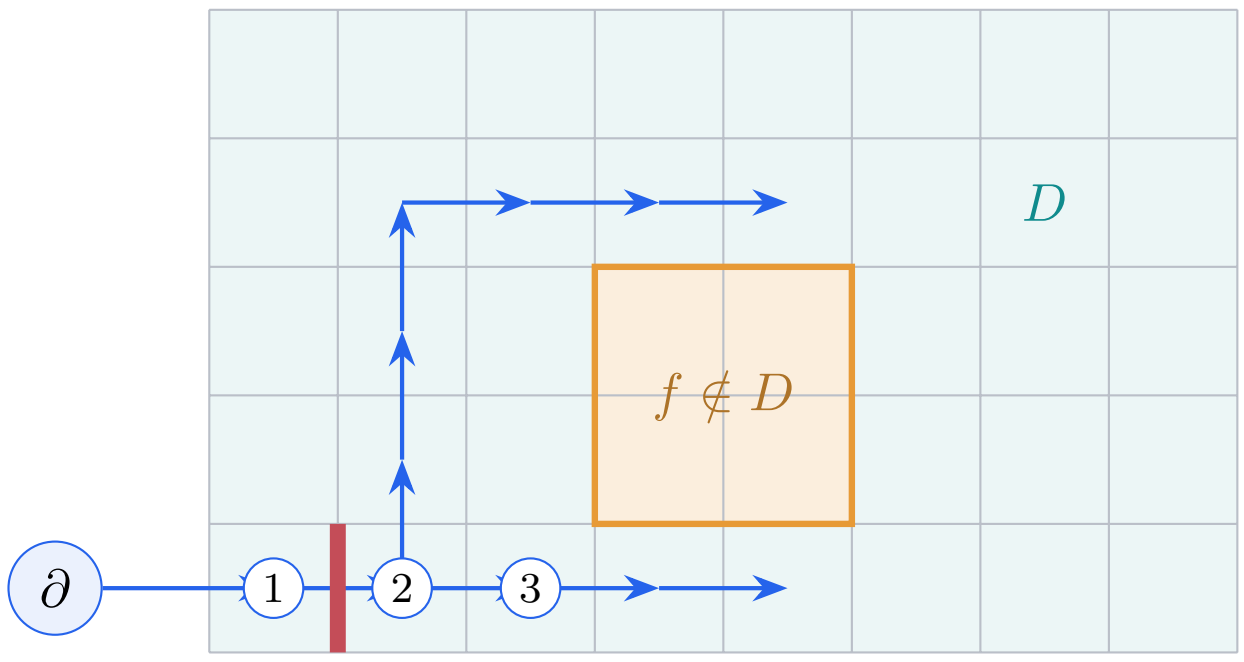

We keep the dual edge between neighboring plaquettes if their common primal edge lies outside $\Gamma$; otherwise we erase it.
(Each such  boundary edge not in $\Gamma$ also gives a dual edge connected to $\partial$.)
Heuristically, the connected components in this dual graph
are simply the unions of plaquettes ``separated by the barriers $\Gamma$''.

Let $D$ be the set of plaquettes reachable from $\partial$. Choose a spanning tree of
$D\cup\{\partial\}$ rooted at $\partial$. 

Now we start to integrate. Process the plaquettes in increasing order of tree
distance from the root $\partial$. At each plaquette $f$, integrate over the primal edge $\tilde e$ crossed by
the dual edge from $f$ to its parent.

The parent is either $\partial$ or a face whose weight has already been integrated.
The chosen edge lies outside $\Gamma$, so it occurs in neither observable $w_p$, $w_{p'}$ nor in any
$H_r$. The {\it only} remaining horizontal weight that uses  $\tilde e$ 
is the weight of $f$, namely 
\[
e^{N \beta_f\Re \Tr(Q^1_f)}e^{N \beta_f\Re \Tr(Q^2_f)}
=e^{t_f (X^1_f+X^2_f)}.
\]
Its variable $g$ appears exactly once in $Q_f$. By left and right invariance of Haar
measure and invariance under inversion, for any $A,B\in G$,
\begin{equ}[e:integrate-const]
 \int_G e^{N\beta_f\operatorname{Re}\Tr(Ag^{\pm1}B)}\,\dH(g)
 =\int_G e^{N\beta_f\operatorname{Re}\Tr u}\,\dH(u).
\end{equ}
This is a constant independent of $A,B$.
Integrating the two replicas separately replaces the  plaquette  weight by the square of
this constant. Continuing along the tree removes every weight in $D$.
It is important to note that for the fixed realization of $R$, these constant factors arising from the above Haar integration over horizontal plaquettes do not depend on the remaining edge variables and can be absorbed into the normalizing constant in \eqref{eq:dens}. Therefore, they play no further role in the  exchange argument below and need not be controlled uniformly
 in $\beta_{12}$.

\medskip

{\bf Step 2.}
For $f\notin D$, let $K$ be its dual component, viewed as a union of plaquettes (each plaquette is viewed as a closed set).
Every edge on the boundary of $K$, including the boundaries of holes and any portion
on the boundary of the whole layer rectangle, belongs to $\Gamma$.  
(Indeed, if $e$ is on the boundary of $K$ and $e$ is not 
on the boundary of the whole layer rectangle, 
then clearly $e$ belongs to $\Gamma$ by definition of dual components; if $e$ is on the boundary of $K$ and is also
on the boundary of the whole layer rectangle, then we also have $e\in \Gamma$ since $f\notin D$.)

The boundary $\partial K$
has even degree at each vertex.
Any such lattice curve, and thus  $\partial K$, decomposes into edge-disjoint simple lattice
circuits.
One of these simple circuits $\gamma\subseteq \Gamma$ in the decomposition therefore encloses $f$.
(This is obvious. To argue more precisely, draw a ray from $f$ along the coordinate direction to the right, and this ray crosses $\partial K$ an odd number of times. So among the simple lattice
circuits, there is a circuit, such that the ray crosses it an odd number of times. Only if the ray goes from the interior of the circuit to the exterior of the circuit, the number of crosses is odd. So this circuit encloses $f$.)

\begin{figure}[h]
\centering
\begin{tikzpicture}[x=.76cm,y=.76cm,>=Latex,every node/.style={font=\small,text=ink},
  note/.style={font=\fontsize{8.5}{11}\selectfont,align=center,fill=white,fill opacity=.90,text opacity=1,inner sep=2pt},scale=0.55]
  \fill[outsidefill] (0,0) rectangle (22,13);
  \fill[insidefill] (2,3) rectangle ++(1,1);
  \fill[insidefill] (2,4) rectangle ++(1,1);
  \fill[insidefill] (2,5) rectangle ++(1,1);
  \fill[insidefill] (2,6) rectangle ++(1,1);
  \fill[insidefill] (2,7) rectangle ++(1,1);
  \fill[insidefill] (2,8) rectangle ++(1,1);
  \fill[insidefill] (2,9) rectangle ++(1,1);
  \fill[insidefill] (2,10) rectangle ++(1,1);
  \fill[insidefill] (3,3) rectangle ++(1,1);
  \fill[insidefill] (3,4) rectangle ++(1,1);
  \fill[insidefill] (3,5) rectangle ++(1,1);
  \fill[insidefill] (3,6) rectangle ++(1,1);
  \fill[insidefill] (3,7) rectangle ++(1,1);
  \fill[insidefill] (3,8) rectangle ++(1,1);
  \fill[insidefill] (3,9) rectangle ++(1,1);
  \fill[insidefill] (3,10) rectangle ++(1,1);
  \fill[insidefill] (4,3) rectangle ++(1,1);
  \fill[insidefill] (4,6) rectangle ++(1,1);
  \fill[insidefill] (4,7) rectangle ++(1,1);
  \fill[insidefill] (4,8) rectangle ++(1,1);
  \fill[insidefill] (4,9) rectangle ++(1,1);
  \fill[insidefill] (4,10) rectangle ++(1,1);
  \fill[insidefill] (5,3) rectangle ++(1,1);
  \fill[insidefill] (5,6) rectangle ++(1,1);
  \fill[insidefill] (5,7) rectangle ++(1,1);
  \fill[insidefill] (5,8) rectangle ++(1,1);
  \fill[insidefill] (5,9) rectangle ++(1,1);
  \fill[insidefill] (5,10) rectangle ++(1,1);
  \fill[insidefill] (6,3) rectangle ++(1,1);
  \fill[insidefill] (6,4) rectangle ++(1,1);
  \fill[insidefill] (6,5) rectangle ++(1,1);
  \fill[insidefill] (6,6) rectangle ++(1,1);
  \fill[insidefill] (6,7) rectangle ++(1,1);
  \fill[insidefill] (6,8) rectangle ++(1,1);
  \fill[insidefill] (6,9) rectangle ++(1,1);
  \fill[insidefill] (6,10) rectangle ++(1,1);
  \fill[insidefill] (7,3) rectangle ++(1,1);
  \fill[insidefill] (7,4) rectangle ++(1,1);
  \fill[insidefill] (7,5) rectangle ++(1,1);
  \fill[insidefill] (7,6) rectangle ++(1,1);
  \fill[insidefill] (7,7) rectangle ++(1,1);
  \fill[insidefill] (7,8) rectangle ++(1,1);
  \fill[insidefill] (7,9) rectangle ++(1,1);
  \fill[insidefill] (7,10) rectangle ++(1,1);
  \fill[insidefill] (8,3) rectangle ++(1,1);
  \fill[insidefill] (8,4) rectangle ++(1,1);
  \fill[insidefill] (8,5) rectangle ++(1,1);
  \fill[insidefill] (8,6) rectangle ++(1,1);
  \fill[insidefill] (8,7) rectangle ++(1,1);
  \fill[insidefill] (9,3) rectangle ++(1,1);
  \fill[insidefill] (9,4) rectangle ++(1,1);
  \fill[insidefill] (9,5) rectangle ++(1,1);
  \fill[insidefill] (9,6) rectangle ++(1,1);
  \fill[insidefill] (9,7) rectangle ++(1,1);
  \fill[otherfill] (4,4) rectangle ++(1,1);
  \fill[otherfill] (4,5) rectangle ++(1,1);
  \fill[otherfill] (5,4) rectangle ++(1,1);
  \fill[otherfill] (5,5) rectangle ++(1,1);
  \fill[otherfill] (13,10) rectangle ++(1,1);
  \fill[otherfill] (13,11) rectangle ++(1,1);
  \fill[otherfill] (14,10) rectangle ++(1,1);
  \fill[otherfill] (14,11) rectangle ++(1,1);
  \fill[otherfill] (15,10) rectangle ++(1,1);
  \fill[otherfill] (15,11) rectangle ++(1,1);
  \fill[otherfill] (16,8) rectangle ++(1,1);
  \fill[otherfill] (16,9) rectangle ++(1,1);
  \fill[otherfill] (17,8) rectangle ++(1,1);
  \fill[otherfill] (17,9) rectangle ++(1,1);
  \fill[otherfill] (18,8) rectangle ++(1,1);
  \fill[otherfill] (18,9) rectangle ++(1,1);
  \draw[grid,line width=.25pt,step=1] (0,0) grid (22,13);
  \draw[muted!70,line width=.5pt] (0,0) rectangle (22,13);
  \draw[wall,line width=1.8pt,line cap=round] (1,6)--(1,7);
  \draw[wall,line width=1.8pt,line cap=round] (1,6)--(2,6);
  \draw[wall,line width=1.8pt,line cap=round] (1,7)--(1,8);
  \draw[wall,line width=1.8pt,line cap=round] (2,8)--(3,8);
  \draw[wall,line width=1.8pt,line cap=round] (3,8)--(4,8);
  \draw[wall,line width=1.8pt,line cap=round] (4,4)--(4,5);
  \draw[wall,line width=1.8pt,line cap=round] (4,4)--(5,4);
  \draw[wall,line width=1.8pt,line cap=round] (4,5)--(4,6);
  \draw[wall,line width=1.8pt,line cap=round] (4,6)--(5,6);
  \draw[wall,line width=1.8pt,line cap=round] (4,7)--(4,8);
  \draw[wall,line width=1.8pt,line cap=round] (5,1)--(5,2);
  \draw[wall,line width=1.8pt,line cap=round] (5,1)--(6,1);
  \draw[wall,line width=1.8pt,line cap=round] (5,2)--(5,3);
  \draw[wall,line width=1.8pt,line cap=round] (5,4)--(6,4);
  \draw[wall,line width=1.8pt,line cap=round] (5,6)--(6,6);
  \draw[wall,line width=1.8pt,line cap=round] (5,11)--(5,12);
  \draw[wall,line width=1.8pt,line cap=round] (5,12)--(6,12);
  \draw[wall,line width=1.8pt,line cap=round] (6,1)--(7,1);
  \draw[wall,line width=1.8pt,line cap=round] (6,4)--(6,5);
  \draw[wall,line width=1.8pt,line cap=round] (6,5)--(6,6);
  \draw[wall,line width=1.8pt,line cap=round] (6,12)--(7,12);
  \draw[wall,line width=1.8pt,line cap=round] (9,1)--(10,1);
  \draw[wall,line width=1.8pt,line cap=round] (10,1)--(11,1);
  \draw[wall,line width=1.8pt,line cap=round] (10,5)--(11,5);
  \draw[wall,line width=1.8pt,line cap=round] (11,1)--(11,2);
  \draw[wall,line width=1.8pt,line cap=round] (11,2)--(12,2);
  \draw[wall,line width=1.8pt,line cap=round] (11,5)--(12,5);
  \draw[wall,line width=1.8pt,line cap=round] (12,2)--(13,2);
  \draw[wall,line width=1.8pt,line cap=round] (12,5)--(12,6);
  \draw[wall,line width=1.8pt,line cap=round] (12,6)--(12,7);
  \draw[wall,line width=1.8pt,line cap=round] (13,2)--(14,2);
  \draw[wall,line width=1.8pt,line cap=round] (13,10)--(13,11);
  \draw[wall,line width=1.8pt,line cap=round] (13,10)--(14,10);
  \draw[wall,line width=1.8pt,line cap=round] (13,11)--(13,12);
  \draw[wall,line width=1.8pt,line cap=round] (13,12)--(14,12);
  \draw[wall,line width=1.8pt,line cap=round] (14,5)--(14,6);
  \draw[wall,line width=1.8pt,line cap=round] (14,5)--(15,5);
  \draw[wall,line width=1.8pt,line cap=round] (14,6)--(14,7);
  \draw[wall,line width=1.8pt,line cap=round] (14,10)--(15,10);
  \draw[wall,line width=1.8pt,line cap=round] (14,12)--(15,12);
  \draw[wall,line width=1.8pt,line cap=round] (15,5)--(16,5);
  \draw[wall,line width=1.8pt,line cap=round] (15,10)--(16,10);
  \draw[wall,line width=1.8pt,line cap=round] (15,12)--(16,12);
  \draw[wall,line width=1.8pt,line cap=round] (16,5)--(17,5);
  \draw[wall,line width=1.8pt,line cap=round] (16,8)--(16,9);
  \draw[wall,line width=1.8pt,line cap=round] (16,8)--(17,8);
  \draw[wall,line width=1.8pt,line cap=round] (16,9)--(16,10);
  \draw[wall,line width=1.8pt,line cap=round] (16,10)--(16,11);
  \draw[wall,line width=1.8pt,line cap=round] (16,10)--(17,10);
  \draw[wall,line width=1.8pt,line cap=round] (16,11)--(16,12);
  \draw[wall,line width=1.8pt,line cap=round] (17,3)--(17,4);
  \draw[wall,line width=1.8pt,line cap=round] (17,4)--(17,5);
  \draw[wall,line width=1.8pt,line cap=round] (17,4)--(18,4);
  \draw[wall,line width=1.8pt,line cap=round] (17,5)--(17,6);
  \draw[wall,line width=1.8pt,line cap=round] (17,5)--(18,5);
  \draw[wall,line width=1.8pt,line cap=round] (17,6)--(17,7);
  \draw[wall,line width=1.8pt,line cap=round] (17,7)--(18,7);
  \draw[wall,line width=1.8pt,line cap=round] (17,8)--(18,8);
  \draw[wall,line width=1.8pt,line cap=round] (17,10)--(18,10);
  \draw[wall,line width=1.8pt,line cap=round] (18,4)--(19,4);
  \draw[wall,line width=1.8pt,line cap=round] (18,5)--(19,5);
  \draw[wall,line width=1.8pt,line cap=round] (18,8)--(19,8);
  \draw[wall,line width=1.8pt,line cap=round] (18,10)--(19,10);
  \draw[wall,line width=1.8pt,line cap=round] (19,5)--(20,5);
  \draw[wall,line width=1.8pt,line cap=round] (19,8)--(19,9);
  \draw[wall,line width=1.8pt,line cap=round] (19,9)--(19,10);
  \draw[wall,line width=1.8pt,line cap=round] (19,9)--(20,9);
  \draw[wall,line width=1.8pt,line cap=round] (20,2)--(21,2);
  \draw[wall,line width=1.8pt,line cap=round] (20,5)--(20,6);
  \draw[wall,line width=1.8pt,line cap=round] (20,6)--(20,7);
  \draw[wall,line width=1.8pt,line cap=round] (20,9)--(21,9);
  \draw[wall,line width=1.8pt,line cap=round] (21,9)--(21,10);
  \draw[wall,line width=1.8pt,line cap=round] (21,10)--(21,11);
  \draw[cycle,line width=2.2pt,line cap=round] (2,3)--(2,4);
  \draw[cycle,line width=2.2pt,line cap=round] (2,3)--(3,3);
  \draw[cycle,line width=2.2pt,line cap=round] (2,4)--(2,5);
  \draw[cycle,line width=2.2pt,line cap=round] (2,5)--(2,6);
  \draw[cycle,line width=2.2pt,line cap=round] (2,6)--(2,7);
  \draw[cycle,line width=2.2pt,line cap=round] (2,7)--(2,8);
  \draw[cycle,line width=2.2pt,line cap=round] (2,8)--(2,9);
  \draw[cycle,line width=2.2pt,line cap=round] (2,9)--(2,10);
  \draw[cycle,line width=2.2pt,line cap=round] (2,10)--(2,11);
  \draw[cycle,line width=2.2pt,line cap=round] (2,11)--(3,11);
  \draw[cycle,line width=2.2pt,line cap=round] (3,3)--(4,3);
  \draw[cycle,line width=2.2pt,line cap=round] (3,11)--(4,11);
  \draw[cycle,line width=2.2pt,line cap=round] (4,3)--(5,3);
  \draw[cycle,line width=2.2pt,line cap=round] (4,11)--(5,11);
  \draw[cycle,line width=2.2pt,line cap=round] (5,3)--(6,3);
  \draw[cycle,line width=2.2pt,line cap=round] (5,11)--(6,11);
  \draw[cycle,line width=2.2pt,line cap=round] (6,3)--(7,3);
  \draw[cycle,line width=2.2pt,line cap=round] (6,11)--(7,11);
  \draw[cycle,line width=2.2pt,line cap=round] (7,3)--(8,3);
  \draw[cycle,line width=2.2pt,line cap=round] (7,11)--(8,11);
  \draw[cycle,line width=2.2pt,line cap=round] (8,3)--(9,3);
  \draw[cycle,line width=2.2pt,line cap=round] (8,8)--(8,9);
  \draw[cycle,line width=2.2pt,line cap=round] (8,8)--(9,8);
  \draw[cycle,line width=2.2pt,line cap=round] (8,9)--(8,10);
  \draw[cycle,line width=2.2pt,line cap=round] (8,10)--(8,11);
  \draw[cycle,line width=2.2pt,line cap=round] (9,3)--(10,3);
  \draw[cycle,line width=2.2pt,line cap=round] (9,8)--(10,8);
  \draw[cycle,line width=2.2pt,line cap=round] (10,3)--(10,4);
  \draw[cycle,line width=2.2pt,line cap=round] (10,4)--(10,5);
  \draw[cycle,line width=2.2pt,line cap=round] (10,5)--(10,6);
  \draw[cycle,line width=2.2pt,line cap=round] (10,6)--(10,7);
  \draw[cycle,line width=2.2pt,line cap=round] (10,7)--(10,8);
  \fill[wall] (1,6) circle (1.15pt);
  \fill[wall] (1,7) circle (1.15pt);
  \fill[wall] (1,8) circle (1.15pt);
  \fill[cycle] (2,3) circle (1.15pt);
  \fill[cycle] (2,4) circle (1.15pt);
  \fill[cycle] (2,5) circle (1.15pt);
  \fill[cycle] (2,6) circle (1.15pt);
  \fill[cycle] (2,7) circle (1.15pt);
  \fill[cycle] (2,8) circle (1.15pt);
  \fill[cycle] (2,9) circle (1.15pt);
  \fill[cycle] (2,10) circle (1.15pt);
  \fill[cycle] (2,11) circle (1.15pt);
  \fill[cycle] (3,3) circle (1.15pt);
  \fill[wall] (3,8) circle (1.15pt);
  \fill[cycle] (3,11) circle (1.15pt);
  \fill[cycle] (4,3) circle (1.15pt);
  \fill[wall] (4,4) circle (1.15pt);
  \fill[wall] (4,5) circle (1.15pt);
  \fill[wall] (4,6) circle (1.15pt);
  \fill[wall] (4,7) circle (1.15pt);
  \fill[wall] (4,8) circle (1.15pt);
  \fill[cycle] (4,11) circle (1.15pt);
  \fill[wall] (5,1) circle (1.15pt);
  \fill[wall] (5,2) circle (1.15pt);
  \fill[cycle] (5,3) circle (1.15pt);
  \fill[wall] (5,4) circle (1.15pt);
  \fill[wall] (5,6) circle (1.15pt);
  \fill[cycle] (5,11) circle (1.15pt);
  \fill[wall] (5,12) circle (1.15pt);
  \fill[wall] (6,1) circle (1.15pt);
  \fill[cycle] (6,3) circle (1.15pt);
  \fill[wall] (6,4) circle (1.15pt);
  \fill[wall] (6,5) circle (1.15pt);
  \fill[wall] (6,6) circle (1.15pt);
  \fill[cycle] (6,11) circle (1.15pt);
  \fill[wall] (6,12) circle (1.15pt);
  \fill[wall] (7,1) circle (1.15pt);
  \fill[cycle] (7,3) circle (1.15pt);
  \fill[cycle] (7,11) circle (1.15pt);
  \fill[wall] (7,12) circle (1.15pt);
  \fill[cycle] (8,3) circle (1.15pt);
  \fill[cycle] (8,8) circle (1.15pt);
  \fill[cycle] (8,9) circle (1.15pt);
  \fill[cycle] (8,10) circle (1.15pt);
  \fill[cycle] (8,11) circle (1.15pt);
  \fill[wall] (9,1) circle (1.15pt);
  \fill[cycle] (9,3) circle (1.15pt);
  \fill[cycle] (9,8) circle (1.15pt);
  \fill[wall] (10,1) circle (1.15pt);
  \fill[cycle] (10,3) circle (1.15pt);
  \fill[cycle] (10,4) circle (1.15pt);
  \fill[cycle] (10,5) circle (1.15pt);
  \fill[cycle] (10,6) circle (1.15pt);
  \fill[cycle] (10,7) circle (1.15pt);
  \fill[cycle] (10,8) circle (1.15pt);
  \fill[wall] (11,1) circle (1.15pt);
  \fill[wall] (11,2) circle (1.15pt);
  \fill[wall] (11,5) circle (1.15pt);
  \fill[wall] (12,2) circle (1.15pt);
  \fill[wall] (12,5) circle (1.15pt);
  \fill[wall] (12,6) circle (1.15pt);
  \fill[wall] (12,7) circle (1.15pt);
  \fill[wall] (13,2) circle (1.15pt);
  \fill[wall] (13,10) circle (1.15pt);
  \fill[wall] (13,11) circle (1.15pt);
  \fill[wall] (13,12) circle (1.15pt);
  \fill[wall] (14,2) circle (1.15pt);
  \fill[wall] (14,5) circle (1.15pt);
  \fill[wall] (14,6) circle (1.15pt);
  \fill[wall] (14,7) circle (1.15pt);
  \fill[wall] (14,10) circle (1.15pt);
  \fill[wall] (14,12) circle (1.15pt);
  \fill[wall] (15,5) circle (1.15pt);
  \fill[wall] (15,10) circle (1.15pt);
  \fill[wall] (15,12) circle (1.15pt);
  \fill[wall] (16,5) circle (1.15pt);
  \fill[wall] (16,8) circle (1.15pt);
  \fill[wall] (16,9) circle (1.15pt);
  \fill[wall] (16,10) circle (1.15pt);
  \fill[wall] (16,11) circle (1.15pt);
  \fill[wall] (16,12) circle (1.15pt);
  \fill[wall] (17,3) circle (1.15pt);
  \fill[wall] (17,4) circle (1.15pt);
  \fill[wall] (17,5) circle (1.15pt);
  \fill[wall] (17,6) circle (1.15pt);
  \fill[wall] (17,7) circle (1.15pt);
  \fill[wall] (17,8) circle (1.15pt);
  \fill[wall] (17,10) circle (1.15pt);
  \fill[wall] (18,4) circle (1.15pt);
  \fill[wall] (18,5) circle (1.15pt);
  \fill[wall] (18,7) circle (1.15pt);
  \fill[wall] (18,8) circle (1.15pt);
  \fill[wall] (18,10) circle (1.15pt);
  \fill[wall] (19,4) circle (1.15pt);
  \fill[wall] (19,5) circle (1.15pt);
  \fill[wall] (19,8) circle (1.15pt);
  \fill[wall] (19,9) circle (1.15pt);
  \fill[wall] (19,10) circle (1.15pt);
  \fill[wall] (20,2) circle (1.15pt);
  \fill[wall] (20,5) circle (1.15pt);
  \fill[wall] (20,6) circle (1.15pt);
  \fill[wall] (20,7) circle (1.15pt);
  \fill[wall] (20,9) circle (1.15pt);
  \fill[wall] (21,2) circle (1.15pt);
  \fill[wall] (21,9) circle (1.15pt);
  \fill[wall] (21,10) circle (1.15pt);
  \fill[wall] (21,11) circle (1.15pt);
  
  \fill[teal!25] (6.06,5.06) rectangle (6.94,5.94);
  \draw[teal,line width=1pt] (6.06,5.06) rectangle (6.94,5.94);
  \node[teal,font=\large] at (6.5,5.5) {$f$};
  \node[teal,font=\Large] at (3.15,9.55) {$K$};
  \node[blue,font=\Large] at (11.7,10.3) {$D$};
  \node[blue,font=\large,draw=blue,fill=white,circle,inner sep=2pt] (root) at (-.85,1.5) {$\partial$};
  \draw[blue,line width=1pt,->] (root)--(.5,1.5)--(2.5,1.5)--(2.5,2.5)--(4.5,2.5);
\end{tikzpicture}
\caption{Another picture of a layer. 
$\Gamma$ consists of 
the orange edges as well as the red circuit which is the simple circuit  $\gamma$ in the
proof that encloses $f$. The set $K$ (light green region) does not include the ``hole'' which is the 2 by 2 square. The boundary $\partial K$ does not include the ``dangling'' orange line inside $K$.}
\end{figure}
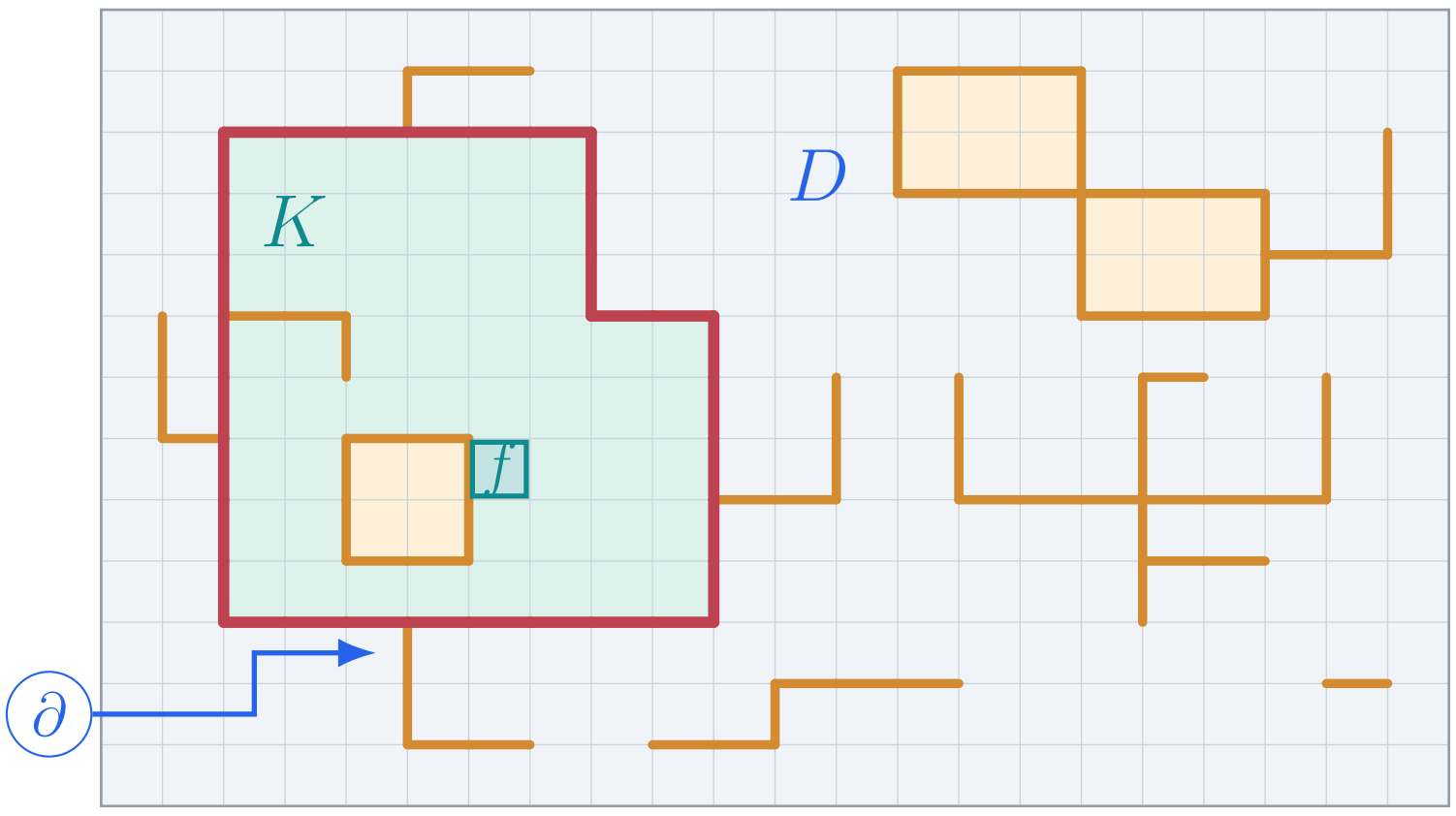

The 3D plaquettes that
supply those edges in $\gamma$ must belong to {\it one} original component $A$, w.r.t. the adjacency relation defined above. 
The circuit
and its interior lie in $B(A)$. In particular, every surviving horizontal weight, that is a horizontal weight factor not integrated in the previous step, has
all $4$ edges in a single original box.

\medskip
{\bf Step 3.}
Recall from \eqref{eq:dens} that, writing $Y=(w_p^1-w_p^2)(w_{p'}^1-w_{p'}^2)$,
\[
\E_{\Prob}[Y\mid R]
=
\frac{\int Y(\mathbf Q)\,L(\mathbf Q)\prod_{r\in R}H_r(\mathbf Q)\,\nu(d \mathbf Q)}
     {\int L(\mathbf Q)\prod_{r\in R}H_r(\mathbf Q)\,\nu(d \mathbf Q)}
\]
Apply the integrations as in the previous steps, in every horizontal layer. Note that the integrated edge variables only appear in $L(\mathbf Q)$, not in $Y(\mathbf Q)$ and $H_r(\mathbf Q)$.
The constants \eqref{e:integrate-const} resulting from these integrations cancel in the numerator and the denominator in the above expression.
Each surviving horizontal weight, each $H_r$
and each observable now has its entire edge support in some $B(A)$. Distinct connected
components of $\mathcal G$ have disjoint unions of boxes. 

If $p\not\leftrightarrow p'$,
exchange $Q_e^1$ and $Q_e^2$ on every remaining edge lying entirely inside the union
of boxes associated with the component of $\mathcal G$ containing $p$.
The product Haar measure and all remaining density factors are invariant. The
integrand in \eqref{eq:rep} changes sign for $F=w_p$, $G=w_{p'}$, so its conditional
expectation 
\[
 \E_{\Prob}\!\left[(F^1-F^2)\overline{(G^1-G^2)}
 \mid R\right] =0,
 \qquad
\mbox{on } p\not\leftrightarrow p'.
\]
On $p\leftrightarrow p'$ its absolute value is at most $4$.
This proves \eqref{eq:switch}. The same exchange works for the conjugated traces
and real parts in the theorem.
\end{proof}

\section{Counting chains of  boxes}
\label{sec:counting}

Suppose $p\ne p'$. On the connection event $p\leftrightarrow p'$
 in \eqref{eq:switch}, choose a simple path $A_1,\ldots,A_k$ in $\mathcal G$ from the component containing $p$ to the one containing $p'$. Recall that this means $B(A_i)\cap B(A_{i+1})\ne\varnothing$. 
 Write $s_i=|A_i|$ and $n=\sum_i s_i$. If both plaquettes lie in the same original component, take $k=1$. The $A_i$ have pairwise disjoint plaquette sets. Of their $n$ plaquettes, all except possibly $p$ and $p'$ must be vertical plaquettes in $R$. Thus, by \eqref{eq:dom}, 
 for any fixed sequence of plaquette sets $A_1,\ldots,A_k$ of total size $n$, the probability that this sequence forms such a path is at most $q^{n-2}$.
 
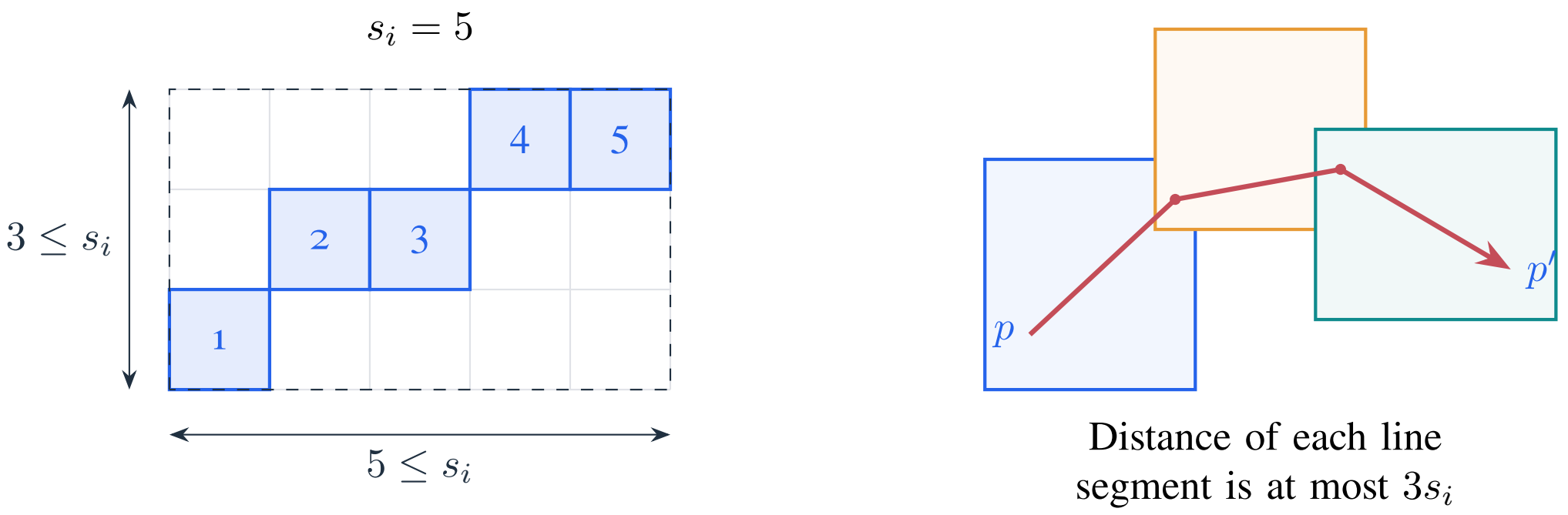
\begin{figure}[h]
\centering
\begin{tikzpicture}[x=.86cm,y=.86cm,>=Stealth,font=\small]
  \begin{scope}
    \draw[muted!20,step=1] (0,0) grid (5,3);
    \foreach \x/\y/\n in {0/0/1,1/1/2,2/1/3,3/2/4,4/2/5}{
      \filldraw[fill=blue!12,draw=blue,thick] (\x,\y) rectangle ++(1,1);
      \node[blue] at (\x+.5,\y+.5) {\n};}
    \draw[ink,dashed] (0,0) rectangle (5,3);
    \draw[<->,ink] (0,-.45)--(5,-.45) node[midway,below] {$5\le s_i$};
    \draw[<->,ink] (-.4,0)--(-.4,3) node[midway,left] {$3\le s_i$};
    \node[align=center] at (2.5,3.6) {$s_i=5$};
  \end{scope}
  \begin{scope}[xshift=7.0cm]
    \draw[blue,thick,fill=blue!6] (0,0) rectangle (2.1,2.3);
    \draw[orange,thick,fill=orange!6] (1.7,1.6) rectangle (3.8,3.6);
    \draw[teal,thick,fill=teal!6] (3.3,.7) rectangle (5.7,2.6);
    \draw[red,very thick,->] (.45,.55)--(1.9,1.9)--(3.55,2.2)--(5.25,1.2);
    \fill[red] (1.9,1.9) circle (.055);
    \fill[red] (3.55,2.2) circle (.055);
    \node[blue,left] at (.45,.55) {$p$};
    \node[blue,right] at (5.25,1.2) {$p'$};
    \node[align=center,text width=5.4cm] at (2.8,-.75)
      {Distance of each line segment is at most $3s_i$};
  \end{scope}
\end{tikzpicture}
\caption{Illustration of \eqref{eq:d}}
\label{fig:distance}
\end{figure}

 A connected set of $s$ plaquettes spans at most $s$ lattice units in each coordinate direction. Indeed, add the plaquettes in the order of a rooted spanning tree, with each parent preceding its children. Each new plaquette shares a vertex with the preceding union and increases each coordinate span by at most one. Hence $B(A_i)$ has $\ell^\infty$ diameter at most $s_i$. See the left picture in Fig.~\ref{fig:distance}.
 Choosing a point in each successive box intersection
 $B(A_i)\cap B(A_{i+1})$
  and applying the triangle inequality in $\ell^1$ gives
\begin{equation}\label{eq:d}
 d(p,p')
 \le \sum_{i} \mbox{diam}_{\ell^1}B(A_i)
 \le 2\sum_i s_i=2n.
\end{equation}
See the right picture in Fig.~\ref{fig:distance}.

To count chains $(A_1,\ldots,A_k)$, we first count the possible choices for a single component.
For the second inequality, the $\ell^1$-diameter of a box is the sum of its coordinate spans.
A box has three coordinate span directions, so we obtain the second inequality with a constant $3$, but we it actually holds with a constant $2$ because of the following reason:
 the first plaquette contributes $2$ to this sum, and each subsequent plaquette can increase it by at most $2$, since only the two spans parallel to the plaquette can grow, each by at most one.

\begin{lemma}\label{lem:count-A}
Fixing a plaquette $p$, the number of possible connected (under adjacency) sets of $s$ plaquettes which contains $p$ is bounded by $48^{2(s-1)}$.
\end{lemma}

\begin{proof}
A lattice vertex belongs to at most $12$ unit squares, and each square has four vertices. Thus the graph (over plaquettes) in which
 plaquettes are adjacent when they share a vertex has degree at most $48$. 
 
 A connected set of $s$ plaquettes containing a specified plaquette $p$ can be encoded by a depth-first traversal of a spanning tree rooted at that plaquette $p$. The walk (i.e. the depth-first traversal) has length $2(s-1)$, and its visited vertices determine the plaquette set. This gives at most $48^{2(s-1)}$ choices.
\end{proof}

\begin{lemma}\label{lem:count-A-chain}
For a fixed  plaquette $p$ contained in $A_1$, the number of possible chains
$(A_1,\ldots,A_k)$ of total size $n$ is at most $D_0^n$, where 
$ D_0=2\cdot48^2\cdot81^2 $.
\end{lemma}

\begin{proof}
Recall $n=\sum_i s_i$. The number of ways to decompose $n$
into $s_1,s_2,...$ is at most $2^{n-1}$.

Use $p$ as the anchor of the first component  $A_1$;
and a (say the lexicographically first) plaquette as the anchor of each subsequent component $A_2,A_3,$ etc. 

A plaquette is specified by a lattice vertex and one of three coordinate planes. 
When consecutive components $A_i,A_{i+1}$ have sizes $s_i$ and $s_{i+1}$, their boxes $B(A_i)$, $B(A_{i+1})$ intersect, so the initial vertices of their anchors are at $\ell^\infty$ distance at most $s_i+s_{i+1}$. There are therefore at most
\[
 3(2(s_i+s_{i+1})+1)^3
 \le
 81^{s_i+s_{i+1}}
\]
choices for the next anchor. 
Multiplying by the bounds for the rooted connected sets in Lemma~\ref{lem:count-A} and summing over the $2^{n-1}$ ordered compositions of $n$ gives 
\[
2^{n-1}
\prod_{i=1}^k 
48^{2(s_i-1)}
\cdot 
\prod_{i=1}^{k-1} 
81^{s_i+s_{i+1}}
\le D_0^n. 
\]
\end{proof}

The above combinatorial count is only a rough upper bound but enough for our purpose. 
In particular: (1) we did not use the condition that  the different components must have disjoint plaquette sets;
(2) we  did not even  impose the condition that 
the chain of boxes  should connect to $p'$. 

However, in the probability estimate  below, we sum only over admissible chains of distinct components with disjoint plaquette sets, which actually reach the distance $d(p,p')$.

\begin{proof}[of Theorem~\ref{theo:main} part (1)]
Set
\[
 q_*=(2D_0)^{-1},\qquad b_*=-\tfrac14\log(1-q_*).
\]
Under \eqref{eq:small} we have $q\le q_*$. Equations \eqref{eq:switch} and \eqref{eq:d} give
\begin{equ}[e:geo-sum]
 |\Cov(w_p,w_{p'})|
 \le 2\sum_{n\ge\max\{2,\lceil d(p,p')/2\rceil\}}D_0^n q_*^{n-2}
 \le 4q_*^{-2}\exp\!\left(-\tfrac{\log 2}{2}\,d(p,p')\right).
\end{equ}
For $p=p'$ the bound follows trivially from $|w_p|\le1$. Multiplication by $N^2$ proves \eqref{eq:main} with $m=(\log 2)/2$. 
In fact this provides explicit choices of $C,b_*,m$ required 
in Theorem~\ref{theo:main}, though not necessarily optimal.
\end{proof}

\section{Higher dimensions}
\label{sec:Higher-dim}

Recall  the
$\lfloor n/2\rfloor$ pairs $(2a-1,2a)$, $1\le a\le\lfloor n/2\rfloor$.
In these orientations the couplings are arbitrary.
Let $T$ be the set of plaquettes of all other orientations. 

The argument in higher dimensions is similar, so we keep the proof brief; however we remark a few points:
\begin{enumerate}
\item
The  coordinate pairs along which the couplings are unrestricted must be {\it disjoint} for the layer argument to apply. 
In four dimensions, for example, we may allow both 
$\beta_{12}$ and $\beta_{34}$ to be arbitrarily large. 
A $12$-layer and a $34$-layer may meet at a vertex, but their {\it edge sets are disjoint}. Each edge in a paired direction therefore belongs to a unique layer.
\item
A marked plaquette may involve several layers. For example, a marked $13$-plaquette has two $1$-edges in distinct $12$-layers and two $3$-edges in distinct $34$-layers. We ``protect'' (when integrating) each of these edges in its own layer. 
We can then perform the planar dual-tree integration separately in each layer and use the ambient boxes $B(A)\subset\Z^n$ for the replica exchange.
\item
The chain count depends on the dimension $n$, while the diameter bound can be proved uniform in $n$. This gives the same decay rate $m=(\log 2)/2$.
\end{enumerate}

\begin{proof}[of Theorem~\ref{theo:high-dim}]
For each paired orientation $(i,j)$, fixing the other $n-2$ coordinates
gives a rectangular $(i,j)$-layer $L$. Each edge in direction $i$ or $j$
lies in exactly one such layer, and these directions belong to no other
pair. Thus distinct layers have disjoint edge sets, even though layers of
different orientations may meet at a vertex.

Mark only the plaquettes of $T$, as in Section~\ref{sec:setup}. With the same $R,H_r$ and
$t_f=N^2\beta_f$, the calculation in \eqref{eq:dom} gives
\[
  \mathbb P(R\supseteq S)\le\prod_{r\in S}(1-e^{-4t_r}),
  \qquad S\subseteq T .
\]
For any realization of $R$ of positive probability, the conditional density
in \eqref{eq:dens} is proportional, relative to product Haar measure, to
\[
  \Big(\prod_{r\in R}H_r\Big)
  \prod_L\ \prod_{f\subseteq L} e^{t_f(X_f^1+X_f^2)} ,
\]
where $f$ ranges over the plaquettes of $L$. Only the marked factors
$H_r$ can involve edge variables from more than one layer.

\smallskip

We do integration in the paired layers.
Form the vertex-connected components $A$ of $R\cup\{p,p'\}$ and their
coordinate boxes $B(A)\subset\Z^n$ as in Section~\ref{sec:layer}. In each layer
$L$, protect the edges used by $R\cup\{p,p'\}$ and call their set $\Gamma$.
For example, a marked $13$-plaquette in four dimensions protects its two
$1$-edges in their respective $12$-layers and its two $3$-edges in their
respective $34$-layers. If $n$ is odd, edges in the unpaired direction
occur only in marked factors and observables.

Apply the exterior dual-tree integration of Section~\ref{sec:layer} within $L$. The
primal edge crossed on the way from a face $f$ to its parent lies outside
$\Gamma$, so it occurs only in paired weights of its unique layer $L$.
The parent is the exterior or a face already processed. Thus the only
remaining weight using this edge is that of $f$, and the Haar integral
over its two copies gives the constant as in \eqref{e:integrate-const}.
The integrations may therefore be performed in any order of the layers.
Their constant factors cancel in the conditional expectation.

Every surviving face in $L$ is enclosed by a circuit in $\Gamma$, by the
planar circuit argument of Section~\ref{sec:layer}. The plaquettes supplying consecutive
circuit edges share a vertex, so all the circuit edges belong to one
$B(A)$. Its intersection with $L$ is a rectangle containing the circuit
and hence the enclosed face. Each remaining weight and each observable
now has its full edge support in one of the ambient boxes. Write
$p\leftrightarrow p'$ for connection in the graph of intersecting boxes,
as before. On disconnection, exchanging the two replicas on the box
component of $p$ preserves all remaining factors and changes the sign of
the integrand in \eqref{eq:rep}. Hence \eqref{eq:switch} holds in dimension $n$, including the
same variants for conjugates and real parts.

\smallskip

Now we do counting.
The argument above \eqref{eq:d} remains exactly the same, so
a chain of total size $\ell$ joining $p$ to
$p'$ satisfies $d(p,p')\le2\ell$.

A vertex belongs to at most $4\binom n2$ plaquettes. The depth-first count
in Lemma~\ref{lem:count-A} thus gives at most $(16\binom n2)^{2(s-1)}$ connected sets
of size $s$ containing a prescribed plaquette. If consecutive components
have sizes $s,u$, their box intersection puts the next anchor within
$\ell^\infty$-distance $s+u$ of the preceding anchor. This gives at most
\[
\binom n2(2s+2u+1)^n\le \Big(3^n\binom n2 \Big)^{s+u}
\] choices. Multiplying and
summing over the ordered compositions as before gives at most $D_n^\ell$
chains of total size $\ell$ from $p$, where we may take
\[
  D_n=2\big(16\textstyle\binom n2\big)^2
         \big(3^n\textstyle\binom n2\big)^2,\qquad
  q_*=(2D_n)^{-1},\qquad b_*(n)=-\tfrac14\log(1-q_*).
\]
For $p\ne p'$, a prescribed simple chain of distinct components joining
them requires at least $\ell-2$ specified marked plaquettes. Its probability is at
most $q_*^{\ell-2}$. Since $D_nq_*=1/2$, the geometric sum in \eqref{e:geo-sum}
now starting at $\ell\ge\max\{2,\lceil d(p,p')/2\rceil\}$, gives
\eqref{eq:disjoint-pairs-decay} with $C_n=4q_*^{-2}$,
recalling that 
$W_p=Nw_p$. 
In case $p= p'$ we simply use $|w_p|\le1$.
%
%
\end{proof}

\section{Infinite-volume limit and Gibbs uniqueness}\label{sec:ym-dlr}

For a finite set $I\subset E^+_{\Z^n}$ of edges and an exterior configuration
$\xi$, let $\gamma_I^\xi$ fix $Q_e=\xi_e$ off $I$ and give the edges in $I$
the law
\[
 \frac{1}{Z_I^\xi}
 \exp \bigg(N\sum_{p:\,\partial p\cap I\ne\varnothing}
                    \beta_p\Re\Tr Q_p\bigg)
 \prod_{e\in I}\dH(Q_e).
\]
Here $p$ ranges over $P^+_{\Z^n}$, and $\partial p$ consists of its four
edges, identified with their positive orientations. A probability measure
on $G^{E^+_{\Z^n}}$ is a Dobrushin--Lanford--Ruelle (DLR) Gibbs measure
if, for every finite $I$, its conditional law given the exterior
configuration $\xi$ is $\gamma_I^\xi$ almost surely. Typically,
exponential clustering does not automatically imply Gibbs uniqueness; one must control the influence from the boundary conditions, c.f. \cite{Borgs84} in the case of lattice gauge theory with finite abelian group.

\begin{theorem}\label{theo:gibbs-uniqueness}
Let $n\ge3$. Fix $N\ge2$ and couplings 
$\boldsymbol\beta$ 
satisfying the assumptions in Theorem~\ref{theo:main} or Theorem~\ref{theo:high-dim}.
The model \eqref{eq:measure} has a unique DLR Gibbs measure
$\mu$ on $\Z^n$. Along every rectangular exhaustion, the finite-volume
laws converge to $\mu$ in total variation on each fixed finite set of edges,
both for free boundaries and uniformly over fixed exterior configurations.
The exponential decay \eqref{eq:main} for $n=3$ and
\eqref{eq:disjoint-pairs-decay} for $n\ge4$ holds under $\mu$ with the
same constants.
\end{theorem}

\begin{proof}
Let $F$ be a bounded measurable function of finitely many edges. Choose
a finite, vertex-connected set $K$ of plaquettes whose boundaries cover
these edges. Write $\Lambda_j=[-j,j]^n\cap\Z^n$ and take $L$ large enough
that $K\subset\Lambda_{\lfloor L/2\rfloor}$. Put $I=E^+_{\Lambda_{2L}}$.
For independent replicas with laws $\gamma_I^\xi$ and $\gamma_I^{\xi'}$,
we have $\gamma_I^\xi(F)-\gamma_I^{\xi'}(F)=\E(F^1-F^2)$.

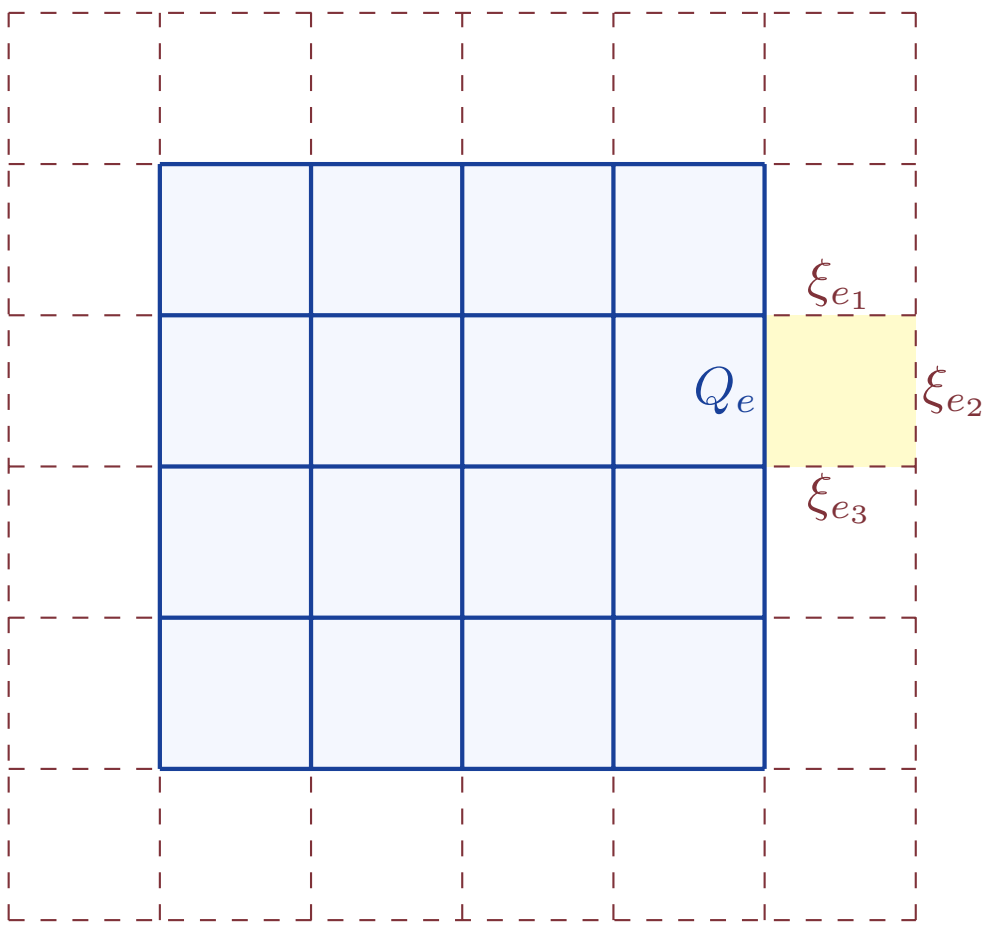
\begin{figure}[htbp]
\centering
\begin{tikzpicture}[x=1cm,y=1cm,font=\small,
  inside/.style={draw=blue!65!black,thick},
  outside/.style={draw=red!65!black,dashed},
  lab/.style={inner sep=1pt}]

  \fill[yellow!25] (2,0) rectangle (3,1);
  \draw[outside] (-3,-3) grid (3,3);
  \fill[blue!5] (-2,-2) rectangle (2,2);
  \draw[inside] (-2,-2) grid (2,2);

  \node[lab,left,text=blue!65!black] at (2,.5) {$Q_e$};

  \node[lab,above,text=red!65!black] at (2.5,1) {$\xi_{e_1}$};
  \node[lab,right,text=red!65!black] at (3,.5)  {$\xi_{e_2}$};
  \node[lab,below,text=red!65!black] at (2.5,0) {$\xi_{e_3}$};
\end{tikzpicture}
\caption{Horizontal slice of $\Lambda_{2L}$ at $x_3=0$ for $L=1$.
The solid blue lines are precisely the edges of $I$ in this layer.
The values on the dashed red edges are fixed by the exterior configuration $\xi$.}
\label{fig:dlr-horizontal-slice}
\end{figure}

Recall that
the ``paired'' orientations are $(2a-1,2a)$, $1\le a\le\lfloor n/2\rfloor$;
below we still call them ``horizontal''.
Let $T$ consist of the plaquettes of $\Lambda_{2L}$ in all other
orientations. (For $n=3$, $T$ consists
of the $13$- and $23$-plaquettes.)
Apply the marking rule of
Section~\ref{sec:setup} to $T$. The 
bound in \eqref{eq:dom} gives $\Prob(S\subseteq R)\le q_*^{|S|}$ for
every $S\subseteq T$ and every $\xi,\xi'$.

Let $\mathcal H=P^+_{\Lambda_{2L}}\setminus T$ be the set of horizontal
plaquettes in the box, and let $\mathcal B$ consist of the plaquettes with
at least one edge in $I$ and at least one edge outside $I$.
Write $\omega_p=e^{t_p(X_p^1+X_p^2)}$. For $p\in\mathcal B$, the exterior
edges in $X_p^1$ are evaluated at $\xi$ and those in $X_p^2$ at $\xi'$.
With $\nu_I(\dd\mathbf Q)=\prod_{e\in I}\dH(Q_e^1)\dH(Q_e^2)$,
the calculation in \eqref{eq:dens} gives
\begin{equation}\label{eq:dlr-first-density}
 \Prob(\dd\mathbf Q\mid R)
 =\frac{1}{Z_R}
  \bigg(\prod_{p\in\mathcal B}\omega_p\bigg)
  \bigg(\prod_{r\in R}H_r\bigg)
  \bigg(\prod_{f\in\mathcal H}\omega_f\bigg)\nu_I(\dd\mathbf Q).
\end{equation}
Here and below the exterior configurations are fixed and suppressed from
the notation. The first product depends on the fields in $I$ only through
edges on the boundary of $\Lambda_{2L}$.


Form the boxes $B(A)$ and the graph $\mathcal G$ from $R\cup K$ as in
Section~\ref{sec:layer}. In each horizontal layer of $\Lambda_{2L}$ that meets
$\Lambda_L$, protect the edges $\Gamma$ used by $R\cup K$. Use the plaquettes
of the layer as dual vertices and join two of them when their common
edge lies outside $\Gamma$. Each such edge belongs to
just one horizontal layer and is not contained in the boundary of
$\Lambda_{2L}$. 

A dual component meeting both the perimeter of its layer
and $\Lambda_L$ spans at least $L$ lattice units in one coordinate and
hence contains at least $L$ plaquettes. See Fig.~\ref{fig:dlr-dual-component}.
Let $D=D(R)$ be the union of the plaquette sets of all such components.
Given $\mathbf Q$ and $R$, independently choose $1$ or $H_f$ in the
expansion \eqref{e:expand-H} at each $f\in D$, with probabilities
proportional to the two summands. Let $M\subseteq D$ be the set of
plaquettes at which $H_f$ is chosen. We call each $f\in D\setminus M$
a constant plaquette. For $f\in D$,
\[
 \Prob(f\notin M\mid\mathbf Q,R)= \frac{1}{1+H_f}
 =e^{-t_f(2+X_f^1+X_f^2)}
 \ge \delta,\qquad
 \delta=\exp\!\left(-4N^2\max_a\beta_{2a-1,2a}\right)>0 .
\]
There are at most $C_nL^n$ dual components in all the layers. Thus the
probability that some component meeting both sets has no constant plaquette
is at most $C_nL^n(1-\delta)^L$.

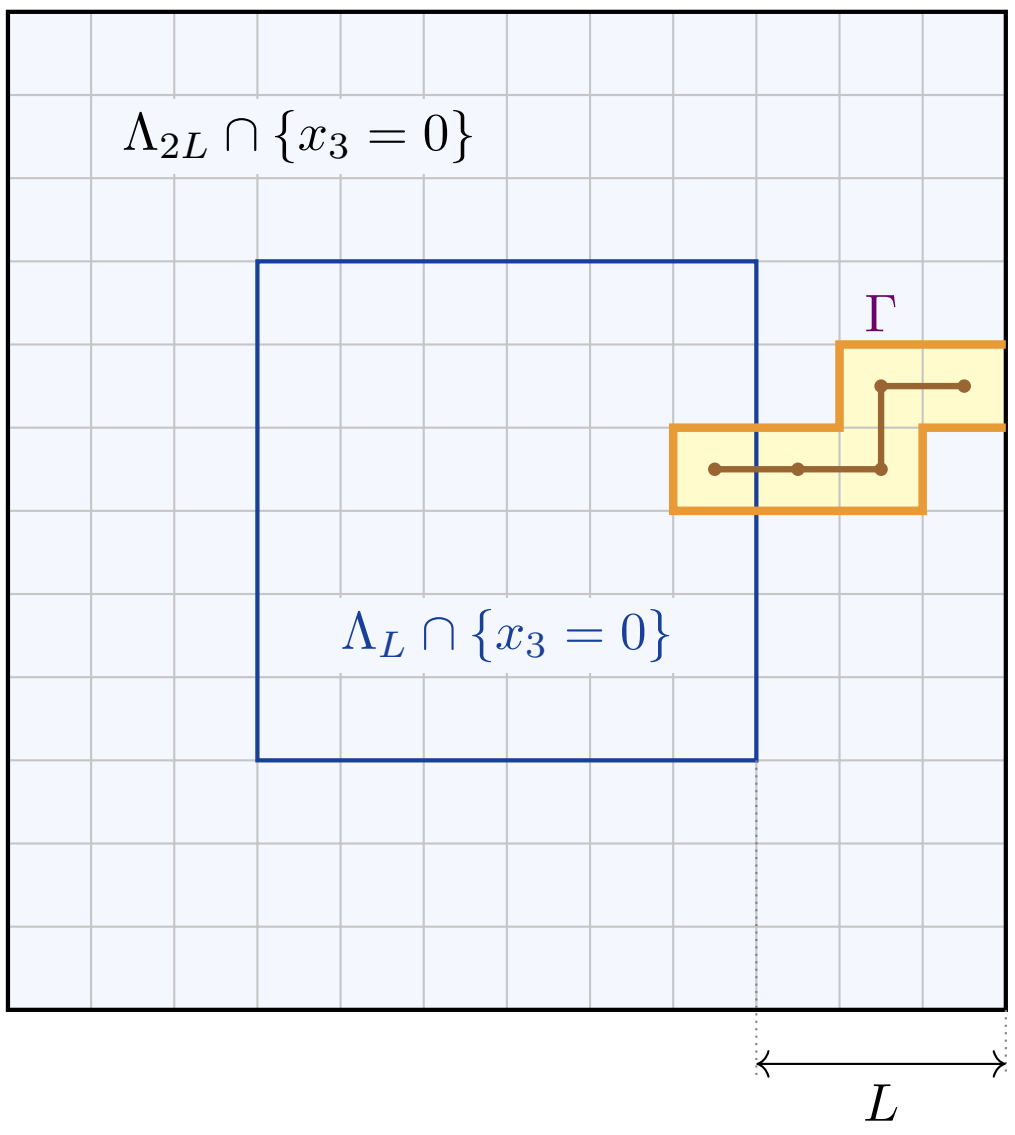
\begin{figure}[htbp]
\centering
\begin{tikzpicture}[scale=.55,font=\small]
  \fill[blue!5] (-6,-6) rectangle (6,6);
  \foreach \x/\y in {2/0,3/0,4/0,4/1,5/1}
    \fill[yellow!25] (\x,\y) rectangle ++(1,1);
  \draw[gray!45,thin] (-6,-6) grid (6,6);
  \draw[black,thick] (-6,-6) rectangle (6,6);
  \draw[blue!65!black,thick] (-3,-3) rectangle (3,3);
  \node[fill=blue!5,inner sep=2pt] at (-2.5,4.5)
    {$\Lambda_{2L}\cap\{x_3=0\}$};
  \node[fill=blue!5,inner sep=2pt,text=blue!65!black] at (0,-1.5)
    {$\Lambda_L \cap\{x_3=0\}$};

  \draw[orange,line width=1.6pt]
    (6,2) -- (4,2) -- (4,1) -- (2,1) -- (2,0)
          -- (5,0) -- (5,1) -- (6,1);
  \node[above,text=violet!85!black,fill=blue!5,inner sep=1.5pt]
    at (4.5,2.05) {$\Gamma$};

  \draw[brown!80!black,very thick]
    (2.5,.5) -- (4.5,.5) -- (4.5,1.5) -- (5.5,1.5);
  \foreach \x/\y in {2.5/.5,3.5/.5,4.5/.5,4.5/1.5,5.5/1.5}
    \fill[brown!80!black] (\x,\y) circle[radius=2.3pt];

  \draw[gray,densely dotted] (3,-3) -- (3,-6.8);
  \draw[gray,densely dotted] (6,-6) -- (6,-6.8);
  \draw[<->] (3,-6.65) -- node[below] {$L$} (6,-6.65);
\end{tikzpicture}
\caption{A horizontal layer for $L=3$. The orange edges belong to
$\Gamma$. 
The yellow plaquettes form a dual component meeting both
$\Lambda_L$ and the right side of the layer perimeter.}
\label{fig:dlr-dual-component}
\end{figure}

Fix a realization of $(R,M)$ of positive probability. In
\eqref{eq:dlr-first-density}, multiplying $\omega_f$ by the probability
of the chosen summand gives $e^{-2t_f}H_f$ for $f\in M$ and
$e^{-2t_f}$ for $f\in D\setminus M$. Absorbing the scalars into
$Z_{R,M}$ gives
\begin{equation}\label{eq:dlr-second-density}
 \Prob(\dd\mathbf Q\mid R,M)
 =\frac{1}{Z_{R,M}}
  \underbrace{\bigg(\prod_{p\in\mathcal B}\omega_p\bigg)
              \bigg(\prod_{r\in R}H_r\bigg)
              \bigg(\prod_{f\in\mathcal H\setminus D}\omega_f\bigg)}
             _{\displaystyle =:\Phi_R(\mathbf Q)}
  \bigg(\prod_{f\in M}H_f\bigg)\nu_I(\dd\mathbf Q).
\end{equation}

Suppose now that each component in $D$ contains a constant plaquette.
In each component root a spanning tree at one such plaquette and process
the other plaquettes in increasing distance from the root.
Each primal edge crossed by a tree edge lies outside $\Gamma$ and off
the boundary of $\Lambda_{2L}$. Its two adjacent horizontal plaquettes
belong to $D$. Hence neither $\Phi_R$ nor $F^1-F^2$ depends on the
two copies of this edge. When a child is processed, its parent's factor
is already constant, so only the child's factor can depend on these copies.
By \eqref{e:integrate-const} and linearity, integrating this factor gives
a scalar for either $1$ or $H_f$. Thus the tree integrations replace
$\prod_{f\in M}H_f$ by the same scalar in the numerator and denominator,
giving
\begin{equation}\label{eq:dlr-tree-expectation}
 \E(F^1-F^2\mid R,M)
 =\frac{\int (F^1-F^2)\Phi_R(\mathbf Q)\,\nu_I(\dd\mathbf Q)}
        {\int \Phi_R(\mathbf Q)\,\nu_I(\dd\mathbf Q)} .
\end{equation}
We may use $\nu_I$ here since the integrands are independent of the
integrated edges and Haar measure has total mass one.

Every remaining horizontal plaquette meeting $\Lambda_L$ belongs to a component
disjoint from the layer perimeter. All boundary edges of this component
are in $\Gamma$. By Section~\ref{sec:layer} a circuit of these edges
encloses the plaquette. Consecutive circuit edges come from vertex-adjacent
plaquettes of $R\cup K$, so the whole plaquette lies in a single $B(A)$.
Each marked factor $H_r$ is already supported in one such box.
If the union of boxes in the component of $\mathcal G$ containing $K$
lies in $\Lambda_L$, exchange the replicas on the remaining edges of this
union. The exchange leaves the boundary factors in $\Phi_R$ untouched.
Every marked or horizontal factor in $\Phi_R$ touched by the exchange has
its full support in the union and is invariant. Since the exchange reverses
$F^1-F^2$, the numerator in \eqref{eq:dlr-tree-expectation} vanishes.
Averaging over the two rounds of marks gives
\[
 \gamma_I^\xi(F)-\gamma_I^{\xi'}(F)
 =\E\big[\E(F^1-F^2\mid R,M)\big].
\]

If this union reaches outside $\Lambda_L$, a simple chain from a fixed
plaquette of $K$ reaches that far. By the box diameter bound in
\eqref{eq:d}, its total size $\ell$ satisfies $2\ell\ge L/2$.
The count in Section~\ref{sec:Higher-dim} gives at most $D_n^\ell$
such chains. For $n=3$ its formula gives $D_3=2\cdot48^2\cdot81^2=D_0$,
the constant in Lemma~\ref{lem:count-A-chain}.
Each prescribed chain requires at least $\ell-|K|$ specified marks, so its
probability under the original joint law is at most $q_*^{\ell-|K|}$.
Since $D_nq_*=1/2$ and $|F^1-F^2|\le 2\|F\|_\infty$, we obtain
\begin{equation}\label{eq:ym-dlr-bound}
 \sup_{\xi,\xi'}|\gamma_I^\xi(F)-\gamma_I^{\xi'}(F)|
 \le 2\|F\|_\infty
 \bigg(C_n L^n(1-\delta)^L+
       q_*^{-|K|}\sum_{\ell\ge L/4}2^{-\ell}\bigg).
\end{equation}
For fixed $N$ and $\boldsymbol\beta$ the right-hand side tends to zero.

In any free box containing $\Lambda_{2L+1}$, conditioning off $I$ writes
the expectation of $F$ as an average of $\gamma_I^\xi(F)$. Thus
\eqref{eq:ym-dlr-bound} bounds the difference of the expectations in any
two such boxes. Taking indicators on a fixed finite set of edges gives a
Cauchy family of marginals in total variation. Their limits are consistent
and define a measure $\mu$ by the Kolmogorov extension theorem. For any
finite edge set $J$ and bounded local $F$, the function
$\xi\mapsto\gamma_J^\xi(F)$ is bounded and local. Hence the identity
\[
 \mu_\Lambda(F)=\int\gamma_J^\xi(F)\,\mu_\Lambda(\dd\xi)
\]
valid in every sufficiently large free box passes to $\mu$. Applying the
identity to products with exterior cylinder indicators gives the DLR
equations by a monotone class argument, since $\gamma_J^\xi$ fixes the
exterior of $J$.

Under any DLR measure, and also in any box containing $\Lambda_{2L+1}$
with a prescribed exterior configuration, the conditional law on $I$
given its exterior $\xi$ is again $\gamma_I^\xi$.
Comparing each law with $\mu$ by \eqref{eq:ym-dlr-bound}, again
using cylinder indicators, proves uniqueness and the stated uniform
convergence. For periodic boxes with all side lengths tending to infinity,
the same conditioning applies because, for each fixed $L$, the cube
$\Lambda_{2L+1}$ eventually embeds in the torus. Passing the corresponding
free-boundary finite volume estimate to $\mu$ gives the claimed exponential decay.
\end{proof}

\section{Staggered vertical couplings in 3D}
\label{sec:ym3-ring}

The proof of Theorem~\ref{theo:3D-ring}
is similar with that of Theorem~\ref{theo:main}, but with an additional trick of integrating ``rings'' to treat the staggered vertical plaquettes with possibly large $\beta_p$. Write convolution on $G$ w.r.t. Haar measure as
\[
 (f*g)(z)=\int_G f(y)g(y^{-1}z)\,dH(y).
 \]

\begin{proof}[of Theorem~\ref{theo:3D-ring}]

As in the proof of Theorem~\ref{theo:main} we again start with 
 two replica of the field $Q^1,Q^2$.

When marking and the vertical plaquettes,
we only apply the marking  {\it only} to the vertical plaquettes
with subscript $b$ or $c$. For their marked set $R$ we have, as before,
\begin{equation}\label{eq:ym3-ring-marks}
 \mathbb P(S\subseteq R)\le q^{|S|},\qquad
 q=1-\exp\big[-4N^2\max\{\beta_b,\beta_c\}\big].
\end{equation}
Fix $R$.  For every vertical plaquette $s$ with $\beta_b$ or $\beta_c$,
\[
\exp\!\left[N\beta_s\bigl(\operatorname{ReTr}Q_s+
                          \operatorname{ReTr}Q'_s\bigr)\right]
=e^{-2N^2\beta_s}(1+H_s).
\]
As before, upon expansion,
the marked plaquettes carry the factors $H_s$, and the unmarked plaquettes with subscript $b$ or $c$ contribute
only factors $1$.

Project the vertical edges of the slab $[k,k+1]$ onto the $12$-plane.:
\begin{itemize}
\item
Each vertical edge from $(i,j,k)$ to $(i,j,k+1)$
 projects to a vertex $(i,j)$;
 \item
each vertical plaquette projects to an
edge joining two such vertices. This edge inherits the ``subscript'' $B,b,C,c$ from the vertical plaquette.
\end{itemize}
The edges with subscript $B$ or $C$ form
vertex-disjoint unit squares with lower left corner $(i,j)$ satisfying
\[
             i\text{ even},\qquad j+k\text{ even}.
\]
Indeed, the $B$-edges pair consecutive columns starting at even $i$, and
the $C$-edges pair consecutive rows starting at $j\equiv k\pmod 2$.

\begin{figure}[htbp]
  \centering
  \begin{tikzpicture}[
    x=.61cm, y=.61cm,
    unrestricted/.style={draw=orange!75!black, line width=1.25pt},
    small coupling/.style={draw=blue!62!black, line width=.55pt,
                           dash pattern=on 2pt off 2pt},
    vertex/.style={fill=black!80},
    font=\small
  ]
    \foreach \parity/\offset/\heading in {0/0/even,1/8/odd} {
      \begin{scope}[shift={(\offset,0)}]
        \node[above] at (2.5,5.23) {$k$ \heading};
        \foreach \i in {0,2,4} {
          \foreach \j in {0,...,4} {
            \ifodd\numexpr\j+\parity\relax
            \else
              \fill[orange!12] (\i,\j) rectangle ++(1,1);
            \fi
          }
        }
        \foreach \j in {0,...,5} {
          \foreach \i in {0,...,4} {
            \ifodd\i\relax
              \draw[small coupling] (\i,\j) -- ++(1,0);
            \else
              \draw[unrestricted] (\i,\j) -- ++(1,0);
            \fi
          }
        }
        \foreach \i in {0,...,5} {
          \foreach \j in {0,...,4} {
            \ifodd\numexpr\j+\parity\relax
              \draw[small coupling] (\i,\j) -- ++(0,1);
            \else
              \draw[unrestricted] (\i,\j) -- ++(0,1);
            \fi
          }
        }
        \foreach \i in {0,...,5}
          \foreach \j in {0,...,5}
            \fill[vertex] (\i,\j) circle[radius=1.05pt];
        \node[below left, font=\scriptsize] at (0,0) {$0$};
        \node[right] at (5.12,0) {$i$};
        \node[left] at (-.10,5) {$j$};
      \end{scope}
    }
  \end{tikzpicture}
  \caption{Projections of the vertical plaquettes in the slab $[k,k+1]$
    onto the $12$-plane. Horizontal edges carry $B$ (solid) or $b$ (dashed).
    Vertical edges in the projection carry $C$ (solid) or $c$ (dashed).
    The shaded unit squares have vertex-disjoint boundaries, each formed by
    two $B$-edges and two $C$-edges. Their positions shift by one row
    when $k$ increases by one.}
  \label{fig:staggered-slab-projection}
\end{figure}
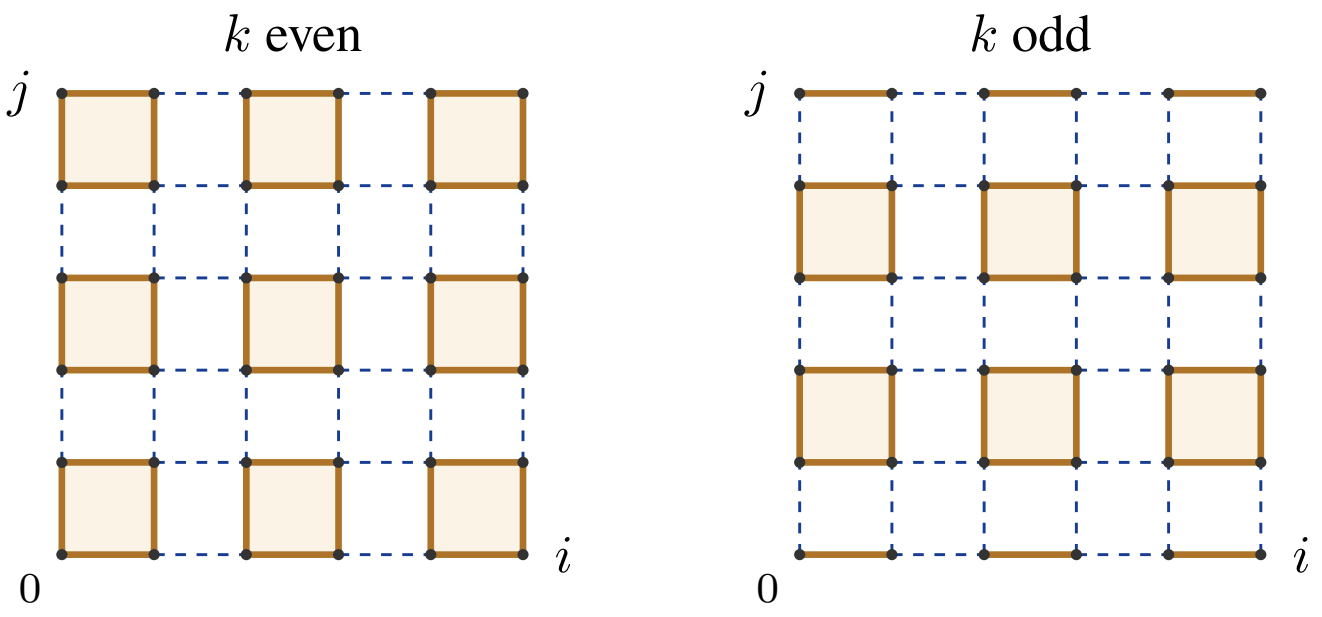

Above each such square, the four side plaquettes of the unit cube form
a {\it ring}, with two $13$ sides and two $23$ sides. 
The four vertical edges of the ring
project to the four corners of the square.
The idea  is to integrate out some of these rings by convolutions, which we describe now.

If any one of the four edges of the ring belongs to a plaquette of $R\cup\{p,p'\}$, keep
all four edges of that ring. Integrate the other rings (separately in each replica),
with the horizontal edges fixed. 
(At a free boundary  of the entire lattice $\Lambda$ a square can be cut down to a single edge or vertex; apply the same retention rule to the remaining edges.)

\smallskip
To describe the convolution on one ring,
define normalized kernel
\[
 k_t(g)=\frac{\exp\{Nt\,\operatorname{Re}\operatorname{Tr}g\}}
              {\int_G\exp\{Nt\,\operatorname{Re}\operatorname{Tr}u\}\,dH(u)},
\]
The normalizing constants are absorbed into the conditional partition function.
Each $k_t$ is central (namely $k_t(g)=k_t(hgh^{-1})$)
 and invariant under inversion.
Number the four corners $0,1,2,3$ around the square starting along a
$13$ side. Let $a_r,b_r$ be the bottom and top edges directed from
corner $r$ to corner $r+1$, and let $v_r$ be the upward edge at corner
$r$, with $v_4=v_0$. The strip shows the ring cut open along $v_0=h$. See Fig.~\ref{fig:staggered-ring-unfolding}.

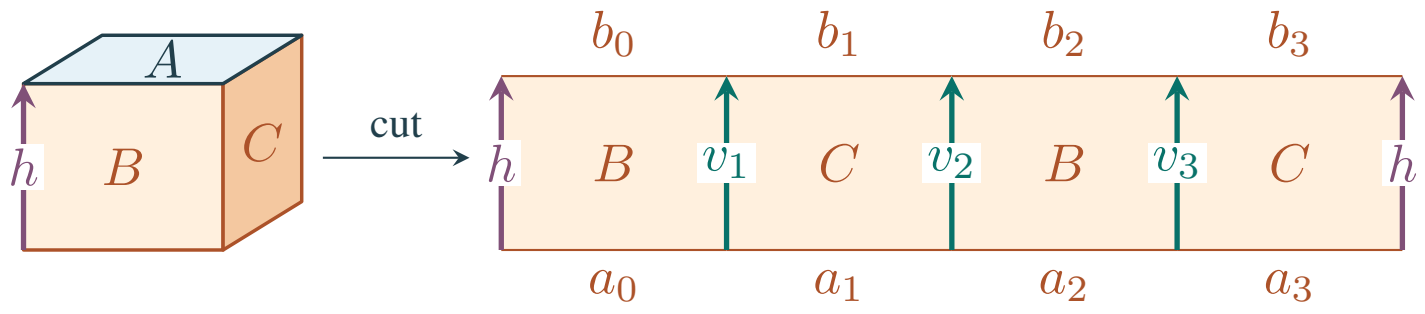
\begin{figure}[htbp]
\centering
\begin{tikzpicture}[x=1cm,y=1cm,font=\small,>=stealth,
    edge label/.style={midway,fill=white,inner sep=.8pt}]
  \definecolor{ringOutline}{HTML}{AC5229}
  \definecolor{ringFront}{HTML}{FFF0DE}
  \definecolor{ringSide}{HTML}{F4C8A0}
  \definecolor{ringTop}{HTML}{E6F2F8}
  \definecolor{ringInk}{HTML}{213F4B}
  \definecolor{ringVertical}{HTML}{087269}
  \definecolor{ringCut}{HTML}{805078}

  \begin{scope}[draw=ringOutline,line join=round,line width=.65pt]
    \filldraw[fill=ringFront] (0,0) rectangle (1.32,1.1);
    \filldraw[fill=ringSide] (1.32,0) -- (1.84,.32)
      -- (1.84,1.42) -- (1.32,1.1) -- cycle;
    \filldraw[draw=ringInk,fill=ringTop] (0,1.1) -- (1.32,1.1)
      -- (1.84,1.42) -- (.52,1.42) -- cycle;
    \node[text=ringOutline] at (.66,.55) {$B$};
    \node[text=ringOutline] at (1.58,.71) {$C$};
    \node[text=ringInk] at (.92,1.26) {$A$};
    \draw[->,draw=ringCut,line width=1pt] (0,0) -- (0,1.1)
      node[edge label,text=ringCut] {$h$};
  \end{scope}

  \draw[->,ringInk] (1.98,.61) -- (2.95,.61)
    node[midway,above,font=\scriptsize] {cut};

  \begin{scope}[shift={(3.16,0)},x=1.49cm,y=1.15cm]
    \foreach \r/\coupling in {0/B,1/C,2/B,3/C} {
      \fill[ringFront] (\r,0) rectangle ++(1,1);
      \draw[ringOutline] (\r,0) -- ++(1,0)
        node[midway,below] {$a_{\r}$};
      \draw[ringOutline] (\r,1) -- ++(1,0)
        node[midway,above] {$b_{\r}$};
      \node[text=ringOutline] at ({\r+.5},.5) {$\coupling$};
    }
    \foreach \r in {1,2,3}
      \draw[->,ringVertical,line width=1pt] (\r,0) -- (\r,1)
        node[edge label,text=ringVertical] {$v_{\r}$};
    \foreach \r in {0,4}
      \draw[->,ringCut,line width=1pt] (\r,0) -- (\r,1)
        node[edge label,text=ringCut] {$h$};
  \end{scope}
\end{tikzpicture}
\caption{Cutting the ring along the vertical edge with variable $h=v_0$.
The four side plaquettes unfold in the order $B,C,B,C$.
The two arrows labeled $h$ represent the same edge.}
\label{fig:staggered-ring-unfolding}
\end{figure}

%
Write $\lambda_0=\lambda_2=\beta_B$, $\lambda_1=\lambda_3=\beta_C$.
The side weights have product
\[
       \prod_{r=0}^3 k_{\lambda_r}
          (a_r v_{r+1} b_r^{-1} v_r^{-1}).
\]
First fix $v_0=h$. We ``bring each $v_r$ back to corner $0$ along
the bottom and top paths'': more precisely we set
\[
 P_r=a_0\cdots a_{r-1},\qquad S_r=b_0\cdots b_{r-1},\qquad
 Y_r=P_r v_r S_r^{-1}\quad(0\le r\le4),
\]
with $P_0=S_0=1$. The lower and upper horizontal holonomies, oriented
the same way, are $U=P_4$ and $V=S_4$. Thus
$Y_0=h$ and $Y_4=UhV^{-1}$. By Haar invariance $dH(v_r)=dH(Y_r)$
for $1\le r\le3$. Moreover,
\[
 Y_{r+1}Y_r^{-1}=P_r(a_r v_{r+1}b_r^{-1}v_r^{-1})P_r^{-1},
 \qquad
 k_{\lambda_r}(a_r v_{r+1}b_r^{-1}v_r^{-1})
     =k_{\lambda_r}(Y_r^{-1}Y_{r+1})
\]
by centrality. 
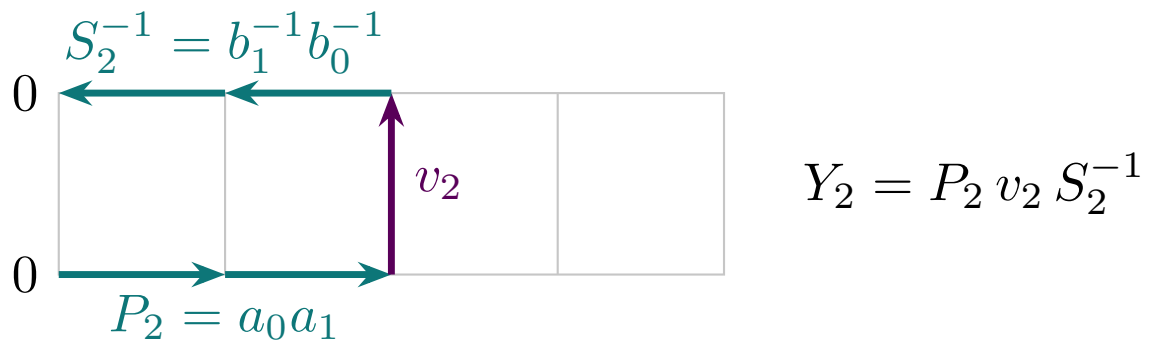
\[
\begin{tikzpicture}[
  x=1.1cm, y=1.2cm, font=\small,
  path arrow/.style={-{Stealth[length=2.2mm]}, line width=1.3pt},
  horizontal/.style={path arrow, teal!85!black},
  vertical/.style={path arrow, violet!70!black}
]
  \draw[gray!45] (0,0) rectangle (4,1);
  \foreach \r in {1,2,3}
    \draw[gray!45] (\r,0) -- (\r,1);

  \draw[horizontal] (0,0) -- (1,0);
  \draw[horizontal] (1,0) -- (2,0);
  \node[below, teal!85!black] at (1,0) {$P_2=a_0a_1$};

  \draw[vertical] (2,0) -- node[right] {$v_2$} (2,1);

  \draw[horizontal] (2,1) -- (1,1);
  \draw[horizontal] (1,1) -- (0,1);
  \node[above, teal!85!black] at (1,1)
    {$S_2^{-1}=b_1^{-1}b_0^{-1}$};

  \node[left] at (0,0) {$0$};
  \node[left] at (0,1) {$0$};
  \node[right] at (4.35,.5) {$Y_2=P_2\,v_2\,S_2^{-1}$};
\end{tikzpicture}
\]
The integration at the first corner is exactly
\[
 \int_G k_{\beta_B}(Y_0^{-1}Y_1)\,
        k_{\beta_C}(Y_1^{-1}Y_2)\,dH(Y_1)
       =(k_{\beta_B}*k_{\beta_C})(Y_0^{-1}Y_2).
\]
Integrating $Y_2,Y_3$ in turn gives $K(h^{-1}UhV^{-1})$, where
$K=k_{\beta_B}*k_{\beta_C}*k_{\beta_B}*k_{\beta_C}$. 
This $K$ is clearly also
is central. Conjugating its argument by $h$ and integrating the fourth
vertical edge gives the factor left by the four sides,
\begin{equation}\label{eq:ym3-ring-kernel}
              T(U,V):=\int_G K(UhV^{-1}h^{-1})\,dH(h).
\end{equation}
The edge variable $h$ transports the top holonomy to the bottom basepoint.
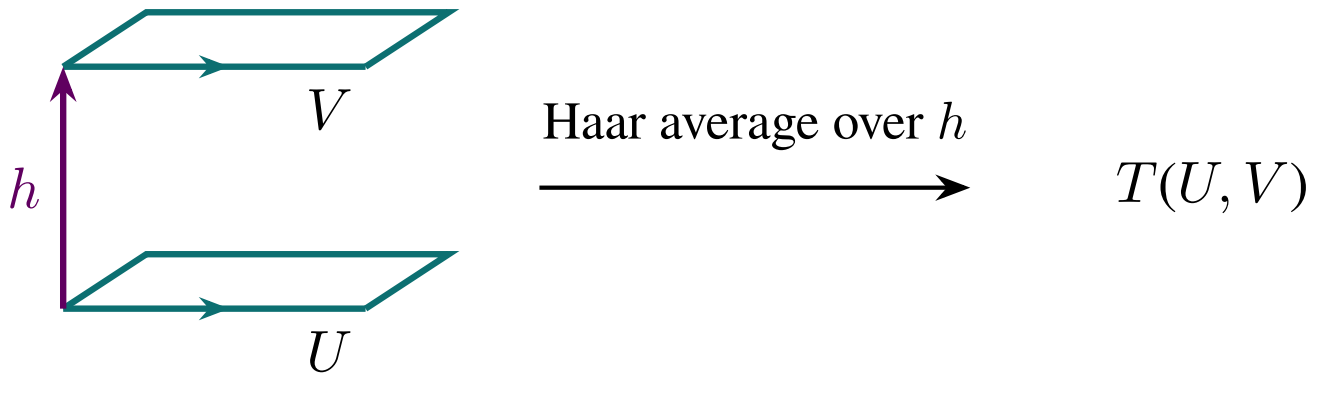
\[
\begin{tikzpicture}[
  loop edge/.style={draw=teal!80!black, line width=1.2pt},
  direction/.style={postaction={decorate}, decoration={markings,
    mark=at position .55 with {\arrow{Stealth[length=2mm]}}}}
]
  \foreach \height/\hol in {0/U,1.6/V} {
    \begin{scope}[yshift=\height cm]
      \draw[loop edge] (0,0) -- (.55,.36) -- (2.55,.36) -- (2,0);
      \draw[loop edge,direction] (0,0) -- (2,0)
        node[pos=.88,below] {$\hol$};
    \end{scope}
  }

  \draw[violet!75!black,line width=1.2pt,-{Stealth[length=2.3mm]}]
    (0,0) -- node[left] {$h$} (0,1.6);

  \draw[-{Stealth[length=2.3mm]},thick] (3.15,.80)
    -- node[above=3pt,font=\small] {Haar average over $h$} (6.0,.80);
  \node at (7.6,.80) {$T(U,V)$};
\end{tikzpicture}
\]
The Haar average makes $T$ invariant under separate conjugations of
$U$ and $V$. Including the bottom and top Wilson weights gives
\begin{equation}\label{eq:ym3-ring-pair}
             F(U,V):=k_{\beta_A}(U)\,k_{\beta_A}(V)\,T(U,V).
\end{equation}
(At the boundary of $\Lambda$, if an unretained projected square consists of a single
edge, integrating one endpoint gives $1$; an isolated vertex also gives $1$.)

\smallskip

Now we have treated the rings, we start to do integrations in each horizontal layer which corresponds to Section~\ref{sec:layer}. The main difference is that now we have ``two-plaquette'' weights $T(U,V)$. To deal with them, we need a notion of ``disjoint pairs''.

\smallskip

Let $f_{i,j,k}$ denote the horizontal plaquette with minimum vertex
$(i,j,k)$. An integrated ring couples $f_{i,j,k}$ to $f_{i,j,k+1}$
only when $i$ is even and $j+k$ is even. A ring immediately below
$f_{i,j,k}$ at the same position would require $j+k-1$ even.
Consequently every horizontal plaquette appears in {\it at most one} pair
\eqref{eq:ym3-ring-pair}. See Fig.~\ref{fig:dis-pair}. The other horizontal plaquettes carry
individual weights $k_{\beta_A}$. This is where the shift of the
$23$ rows with height enters the proof.

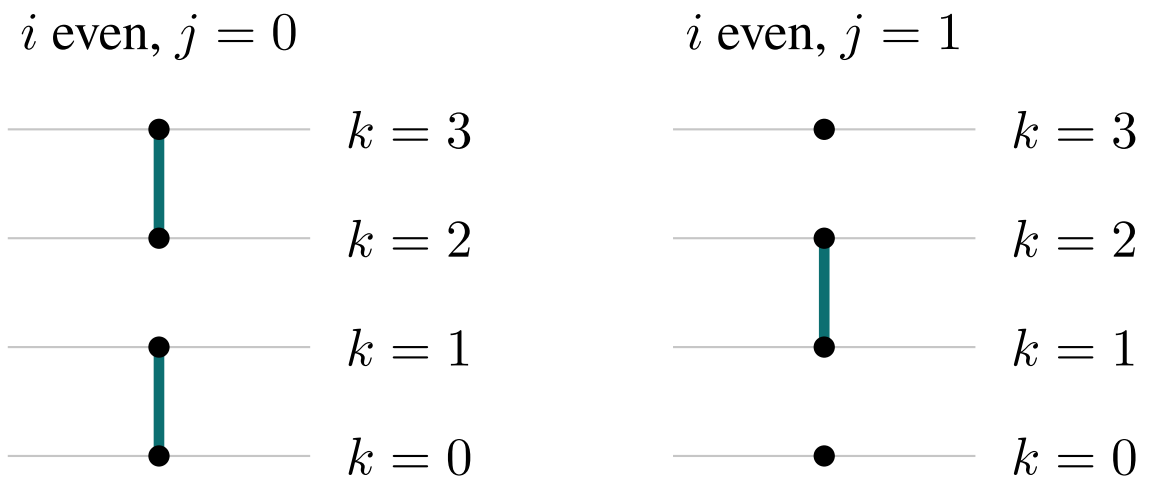
\begin{figure}[h]
\centering
\begin{tikzpicture}[x=1cm,y=.72cm,font=\small]
  \foreach \j/\offset in {0/0,1/4.4} {
    \begin{scope}[xshift=\offset cm]
      \node at (0,3.85) {$i$ even, $j=\j$};

      \foreach \k in {0,...,3} {
        \draw[gray!45,thin] (-1,\k) -- (1,\k);
        \node[right] at (1.1,\k) {$k=\k$};
      }

      \foreach \k in {0,...,2} {
        \ifodd\numexpr\j+\k\relax\else
          \draw[teal!80!black,line width=2pt] (0,\k) -- (0,{\k+1});
        \fi
      }

      \foreach \k in {0,...,3}
        \fill (0,\k) circle[radius=2pt];
    \end{scope}
  }
\end{tikzpicture}
\caption{Illustration for ``every horizontal plaquette appears
in at most one pair''. Here each black dot stands for a horizontal plaquette (not vertex!)}
\label{fig:dis-pair}
\end{figure}

For each $s\in R\cup\{p,p'\}$ let $B_s$ be its coordinate box enlarged
by one lattice unit on each side. A ring retained on account of $s$
lies entirely in $B_s$.
(For instance if $s$ is a plaquette of size $0\times 1\times 1$ then
$B_s$ is of size $2\times 3\times 3$.)
 In each horizontal layer, we protect the edges
contained in an individual $B_s$ and run the exterior-tree integration
of Section~\ref{sec:layer}.  More precisely,
in the horizontal layer of height $k$,
writing $E_k$ for the set of horizontal edges in this layer,
the counterpart of the set $\Gamma$ in Section~\ref{sec:layer} is now defined as
\[
\Gamma_k
=\bigcup_{s\in R\cup\{p,p'\}}\{e\in E_k:e\subseteq B_s\}.
\]
When the spanning tree reaches a paired plaquette, its parent is either $\partial$
or already processed, and its partner is in a different layer (layer $k-1$ or layer $k+1$). The parent
edge $g$ therefore occurs only in the plaquette being processed. In one replica,
\[
 \int_G F(u g^{\pm1}v,V)\,dH(g)=\int_G F(z,V)\,dH(z).
\]
This leaves a weight on the partner alone. 
Integrating that argument again
gives a {\it constant}, so we are again in the situation \eqref{e:integrate-const}.
This ``localization'' (i.e. no third layer gets involved) uses only the disjointness of the pairs, and does not require bound on $T(U,V)$.

Join the plaquettes of $R\cup\{p,p'\}$ whose boxes $B_s$ intersect;
this replaces the notion of adjacency in Section~\ref{sec:layer}.
For each component $J$ let $B(J)$ be the coordinate hull of its boxes.
By the circuit argument of Section~\ref{sec:layer}, every surviving horizontal
plaquette lies in some $B(J)$. If both plaquettes of a pair survive, they are
vertical translates by one unit. Their joint factor is therefore supported
in $B(J)$ enlarged by one further unit, which we denote by $\widehat B(J)$.
The retained side weights, marked factors and observations were already
supported in individual $B_s$.

Connect components whose $\widehat B(J)$ intersect. As before, conditional
on $R$, exchanging replicas on the component of $p$ cancels the covariance
integrand whenever $p$ and $p'$ are disconnected. Hence for
$w_p=N^{-1}W_p$ the covariance is bounded in absolute value by
$2\mathbb P(p\leftrightarrow p')$. To apply the count of Section~\ref{sec:counting}, observe
that the $B_s$ have bounded intersection degree and that a component of
$\ell$ plaquettes has enlarged box of side lengths at most $3\ell+2\le5\ell$.
A simple chain with $n$ plaquettes therefore spans distance at most $15n$.
Consecutive components of sizes $\ell,\ell'$ have anchors within
$5(\ell+\ell')$ in $\ell^\infty$. The same count now gives at most $D^n$
chains from $p$ for an absolute constant $D$. All but at most two plaquettes of an
admissible chain must be marked. Choose $\delta$ so that
$1-e^{-4\delta}\le(2D)^{-1}$. The geometric sum in Section~\ref{sec:counting},
together with \eqref{eq:ym3-ring-marks} and $W_p=Nw_p$, gives
\eqref{eq:ym3-ring-decay}. 
%
\end{proof}

\bibliographystyle{./Martin}
\bibliography{./refs}

\def\cprime{$'$} \def\polhk#1{\setbox0=\hbox{#1}{\ooalign{\hidewidth
  \lower1.5ex\hbox{`}\hidewidth\crcr\unhbox0}}}
\begin{thebibliography}{vdBM94}
\def\myhref#1#2{\href{#2}{\nolinkurl{#1}}}

\bibitem[AC25]{AC25}
\textsc{A.~Adhikari} and \textsc{S.~Cao}.
\newblock Correlation decay for finite lattice gauge theories at weak coupling.
\newblock \emph{The Annals of Probability} \textbf{53}, (2025), 140--174.
\newblock
  \myhref{doi:10.1214/24-AOP1702}{https://dx.doi.org/10.1214/24-AOP1702}.

\bibitem[Aiz82]{Aizenman1982}
\textsc{M.~Aizenman}.
\newblock Geometric analysis of {$\phi^4$} fields and {Ising} models. {Parts I
  and II}.
\newblock \emph{Communications in Mathematical Physics} \textbf{86}, no.~1,
  (1982), 1--48.
\newblock
  \myhref{doi:10.1007/BF01205659}{https://dx.doi.org/10.1007/BF01205659}.

\bibitem[BCD97]{billo1997}
\textsc{M.~Bill{\'o}}, \textsc{M.~Caselle}, and \textsc{A.~D'Adda}.
\newblock Analytic results in (2+1)-dimensional finite temperature {LGT}.
\newblock \emph{International Journal of Modern Physics A} \textbf{12}, no.~32,
  (1997), 5753--5766.

\bibitem[Bor84]{Borgs84}
\textsc{C.~Borgs}.
\newblock Translation symmetry breaking in four dimensional lattice gauge
  theories.
\newblock \emph{Communications in Mathematical Physics} \textbf{96}, no.~2,
  (1984), 251--284.
\newblock
  \myhref{doi:10.1007/BF01240222}{https://dx.doi.org/10.1007/BF01240222}.

\bibitem[CC98]{CampbellChayes1998}
\textsc{M.~Campbell} and \textsc{L.~Chayes}.
\newblock The isotropic {$O(3)$} model and the {Wolff} representation.
\newblock \emph{Journal of Physics A: Mathematical and General} \textbf{31},
  no.~13, (1998), L255--L259.
\newblock
  \myhref{doi:10.1088/0305-4470/31/13/002}{https://dx.doi.org/10.1088/0305-4470/31/13/002}.

\bibitem[Che19]{Chevyrev19YM}
\textsc{I.~Chevyrev}.
\newblock Yang-{M}ills measure on the two-dimensional torus as a random
  distribution.
\newblock \emph{Comm. Math. Phys.} \textbf{372}, no.~3, (2019), 1027--1058.
\newblock
  \myhref{doi:10.1007/s00220-019-03567-5}{https://dx.doi.org/10.1007/s00220-019-03567-5}.

\bibitem[CMR98]{ChayesMachtaRedner1998}
\textsc{L.~Chayes}, \textsc{J.~Machta}, and \textsc{O.~Redner}.
\newblock Graphical representations for {Ising} systems in external fields.
\newblock \emph{Journal of Statistical Physics} \textbf{93}, no. 1--2, (1998),
  17--32.
\newblock
  \myhref{arXiv:cond-mat/9806312}{https://arxiv.org/abs/cond-mat/9806312}.
\newblock
  \myhref{doi:10.1023/B:JOSS.0000026726.43558.80}{https://dx.doi.org/10.1023/B:JOSS.0000026726.43558.80}.

\bibitem[DC18]{DuminilCopin2018}
\textsc{H.~Duminil-Copin}.
\newblock Random currents expansion of the {Ising} model.
\newblock In \emph{European Congress of Mathematics},  869--889. European
  Mathematical Society, 2018.
\newblock \myhref{arXiv:1607.06933}{https://arxiv.org/abs/1607.06933}.
\newblock \myhref{doi:10.4171/176-1/39}{https://dx.doi.org/10.4171/176-1/39}.

\bibitem[Dri89]{Driver89}
\textsc{B.~K. Driver}.
\newblock Y{M{${}_2$}}: continuum expectations, lattice convergence, and
  lassos.
\newblock \emph{Comm. Math. Phys.} \textbf{123}, no.~4, (1989), 575--616.
\newblock
  \myhref{doi:10.1007/BF01218586}{https://dx.doi.org/10.1007/BF01218586}.

\bibitem[ES88]{EdwardsSokal1988}
\textsc{R.~G. Edwards} and \textsc{A.~D. Sokal}.
\newblock Generalization of the {Fortuin-Kasteleyn-Swendsen-Wang}
  representation and {Monte Carlo} algorithm.
\newblock \emph{Physical Review D} \textbf{38}, no.~6, (1988), 2009--2012.
\newblock
  \myhref{doi:10.1103/PhysRevD.38.2009}{https://dx.doi.org/10.1103/PhysRevD.38.2009}.

\bibitem[For22]{For22}
\textsc{M.~P. Forsstr{\"o}m}.
\newblock Decay of correlations in finite {Abelian} lattice gauge theories.
\newblock \emph{Communications in Mathematical Physics} \textbf{393}, (2022),
  1311--1346.
\newblock
  \myhref{doi:10.1007/s00220-022-04391-0}{https://dx.doi.org/10.1007/s00220-022-04391-0}.

\bibitem[GHS70]{GriffithsHurstSherman1970}
\textsc{R.~B. Griffiths}, \textsc{C.~A. Hurst}, and \textsc{S.~Sherman}.
\newblock Concavity of magnetization of an {Ising} ferromagnet in a positive
  external field.
\newblock \emph{Journal of Mathematical Physics} \textbf{11}, no.~3, (1970),
  790--795.
\newblock \myhref{doi:10.1063/1.1665211}{https://dx.doi.org/10.1063/1.1665211}.

\bibitem[Gio15]{giordano2015}
\textsc{M.~Giordano}.
\newblock Renormalization of gauge theories on general anisotropic lattices and
  high-energy scattering in {QCD}.
\newblock \emph{Physical Review D} \textbf{92}, no.~3, (2015), 034514.

\bibitem[Kar82]{karsch1982n}
\textsc{F.~Karsch}.
\newblock {SU(N)} gauge theory couplings on asymmetric lattices.
\newblock \emph{Nuclear Physics B} \textbf{205}, no.~2, (1982), 285--300.

\bibitem[L{\'e}v17]{Levy11}
\textsc{T.~L{\'e}vy}.
\newblock The master field on the plane.
\newblock \emph{Ast\'{e}risque} , no. 388, (2017), ix+201.

\bibitem[Nis25]{Nissim25}
\textsc{R.~Nissim}.
\newblock {$\mathrm{U}(N)$} lattice {Yang--Mills} in the {'t Hooft} regime.
\newblock \emph{arXiv preprint arXiv:2510.22788} (2025).
\newblock
  \myhref{doi:10.48550/arXiv.2510.22788}{https://dx.doi.org/10.48550/arXiv.2510.22788}.

\bibitem[Orl07]{Orland07}
\textsc{P.~Orland}.
\newblock Glueball masses in $(2+1)$-dimensional anisotropic weakly-coupled
  {Yang--Mills} theory.
\newblock \emph{Physical Review D} \textbf{75}, (2007), 101702(R).
\newblock
  \myhref{doi:10.1103/PhysRevD.75.101702}{https://dx.doi.org/10.1103/PhysRevD.75.101702}.

\bibitem[OS78]{OS1978}
\textsc{K.~Osterwalder} and \textsc{E.~Seiler}.
\newblock Gauge field theories on a lattice.
\newblock \emph{Annals of Physics} \textbf{110}, no.~2, (1978), 440--471.

\bibitem[Per75]{Percus1975}
\textsc{J.~K. Percus}.
\newblock Correlation inequalities for {Ising} spin lattices.
\newblock \emph{Communications in Mathematical Physics} \textbf{40}, no.~3,
  (1975), 283--308.
\newblock
  \myhref{doi:10.1007/BF01610004}{https://dx.doi.org/10.1007/BF01610004}.

\bibitem[PS26]{piazza2026log}
\textsc{J.~Piazza} and \textsc{H.~Shen}.
\newblock Log-{S}obolev inequality for {2D $SU(N)$} lattice gauge theory at all
  couplings.
\newblock \emph{arXiv preprint arXiv:2610.35286} (2026).

\bibitem[Sei82]{Seiler}
\textsc{E.~Seiler}.
\newblock \emph{Gauge theories as a problem of constructive quantum field
  theory and statistical mechanics}, vol. 159 of \emph{Lecture Notes in
  Physics}.
\newblock Springer-Verlag, Berlin, 1982,  v+192.

\bibitem[SZZ23]{SZZ22}
\textsc{H.~Shen}, \textsc{R.~Zhu}, and \textsc{X.~Zhu}.
\newblock A stochastic analysis approach to lattice {Y}ang-{M}ills at strong
  coupling.
\newblock \emph{Comm. Math. Phys.} \textbf{400}, no.~2, (2023), 805--851.
\newblock
  \myhref{doi:10.1007/s00220-022-04609-1}{https://dx.doi.org/10.1007/s00220-022-04609-1}.

\bibitem[vdBM94]{vanDenBergMaes1994}
\textsc{J.~van~den Berg} and \textsc{C.~Maes}.
\newblock Disagreement percolation in the study of {Markov} fields.
\newblock \emph{The Annals of Probability} \textbf{22}, no.~2, (1994),
  749--763.
\newblock
  \myhref{doi:10.1214/aop/1176988728}{https://dx.doi.org/10.1214/aop/1176988728}.

\end{thebibliography}

\end{document}